\documentclass{amsart}
\usepackage{amsfonts}
\usepackage{amssymb}
\newtheorem{thm}{Theorem}[section]

\newtheorem{lem}[thm]{Lemma}
\newtheorem{prop}[thm]{Proposition}
\newtheorem{example}[thm]{Example}
\theoremstyle{definition}
\newtheorem{define}[thm]{Definition}
\theoremstyle{remark}

\numberwithin{equation}{section}

\begin{document}

\title[Topological Classification of Integrable Systems]{Symplectic Topology
of Integrable Hamiltonian Systems, II: Topological Classification}

\author{Nguyen Tien Zung}
\address{D\'epartement de Math\'ematiques, Universit\'e Montpellier II}
\email{tienzung@math.univ-montp2.fr}

\date{October 2000, V.2 May 2001}


\dedicatory{Dedicated to Professor Charles-Michel Marle
}
\begin{abstract}

The main purpose of this paper is to give a topological and symplectic
classification of completely integrable Hamiltonian systems in terms of
characteristic classes and other local and global invariants.

\end{abstract}
\maketitle


\section{Introduction}

In this paper we are concerned with topological aspects of completely integrable
Hamiltonian systems, or integrable systems in short. The problem of studying local
and global topological properties of integrable systems is a very natural problem,
with many possible applications, and it has attracted many mathematicians over the
past decades. There are even a few recent books on the subject (see e.g.
\cite{BF,CB,LeUm5,Audin5,Fomenko5}).

The purpose of this paper is to give a topological and symplectic classification of
integrable systems in terms of characteristic classes and other local and global
invariants. Before trying to formulate our theorems, let us recall some of the main
results that have been obtained to date in this direction.

One of the most remarkable results is due to Duistermaat \cite{Duistermaat}, who
defined the {\it monodromy} and the {\it Chern class} of the {\it regular part} of
an integrable system. The monodromy phenomenon has since then been studied by many
authors, for both classical and quantum integrable systems (see e.g.
\cite{CB,CK,GU,Vu1,Vu3}). It is observed in \cite{ZungFocus} that the so-called {\it
focus-focus singularities} are the main (if not the only practical) source of
non-trivial monodromy for ``real-world'' integrable systems. The work of Duistermaat
was developed further and extended to the case of complete isotropic fibrations by
Dazord and Delzant \cite{DD}, and extended by Boucetta and Molino \cite{BM} to
include nondegenerate elliptic singularities of integrable systems.

On the other hand, Fomenko and his collaborators (see e.g. \cite{BF,Fomenko5})
developed a {\it Morse theory} for integrable systems, which takes into account {\it
corank-1} singularities. Fomenko and his school studied mainly systems with two
degrees of freedom, though some results are also valid for higher-dimensional
systems. In particular, they obtained a complete topological classification of
nondegenerate integrable systems on isoenergy 3-manifolds.

The main drawback of both Duistermaat's and Fomenko's theories is the absence of
higher-corank non-elliptic singularities in their picture. It is not surprising,
since the topological structure of higher-corank non-elliptic singularities was
almost completely unknown until more recently. This drawback is quite significant,
because interesting things often happen at singularities, and a lot of global
information is contained there.

Lerman and Umanskii \cite{LeUm5} were among the first people who attacked the
problem of describing the topology of corank-2 nondegenerate singularities of
integrable systems. However, their description is a little bit too complicated in
our view, and they restricted their attention to the case with only one fixed point
in a system with two degrees of freedom.

In a series of papers \cite{ZungAL,ZungFocus,ZungCorank1,ZungNormal}, we studied
local and semi-local aspects of singularities integrable systems. Our main results
include: the existence of a local converging Birkhoff normal form for any analytic
integrable system, the existence of local torus actions and partial action-angle
coordinates, and the topological decomposition of nondegenerate singularities of
higher corank into almost direct products of simplest singularities (elliptic and
hyperbolic corank-1 and focus-focus corank-2). In particular, we gave a much simpler
description of singularities studied by Lerman and Umanskii.

What we do in this paper is to combine our knowledge of singularities with the ideas
of Duistermaat, Fomenko and others in order to study global aspects of ``generic''
integrable systems with singularities. In particular, we will develop the notions of
monodromy and Chern class to take singularities into account. It was a non-trivial
task, because we don't know of a general recipe to define characteristic classes for
singular fibrations: In our case in general there are no ``local sections'' or
``trivial systems'' to speak of, so the obstruction theory does not work directly.
One may try to define some kind of classifying space and universal system, but we
have no idea what they might look like at the moment. And one may try to use the
sheaf of local automorphisms, but this is a very big non-Abelian sheaf, not easy to
deal with. So our first attempts at defining characteristic classes \cite{Zung95}
were not very successful. We then realized that a more detailed study of the sheaf
of local automorphisms of an integrable system allows us to reduce this non-Abelian
structural sheaf to a very nice finite-dimensional Abelian subsheaf, which we will
call the {\it affine monodromy sheaf}. Our main characteristic class, which will be
called the {\it Chern class} and which classifies ``generic'' integrable systems
topologically, is an element of the second cohomology group of this affine monodromy
sheaf. In the regular part of the system, the affine monodromy sheaf is essentially
the same as the monodromy defined by Duistermaat, and our Chern class is also
essentially the same as the one defined by Duistermaat. In the case of
2-degree-of-freedom systems on isoenergy 3-manifolds studied by Fomenko et al., our
affine monodromy contains information about the ``marks'' in the so-called ``marked
molecules'' in Fomenko's classification.

Let us now formulate the main results of this paper.

Let ${\bf F} = (F_1,\dots, F_n) : (M^{2n}, \omega) \to {\Bbb R}^n$ be a smooth
moment map of a completely integrable Hamiltonian system with $n$ degrees of
freedom. We will always assume ${\bf F}$ to be a proper map. In particular, regular
connected components of the level sets of ${\bf F}$ are Lagrangian tori, according
to a classical theorem of Liouville. Denote by $O$ the space of connected components
of the level sets of $F$. We will call $O$ the {\it base space}, or also the {\it
orbit space}, of the system. We have a projection $\pi$ from $(M^{2n},\omega)$ to
$O$ and a map $\tilde{\bf F}: O \to {\mathbb R}^n$, such that ${\bf F} = \tilde{\bf
F} \circ \pi$. The topology of $O$ is induced from $M^{2n}$. For ``generic''
integrable systems, $O$ is a stratified $n$-dimensional manifold. We may view $\pi:
(M^{2n},\omega) \to O$ as the projection map of a singular Lagrangian torus
fibration. We will call it the {\it associated singular Lagrangian fibration} of the
system, and denote it by ${\mathcal L}$: fibers of $\mathcal L$ are connected
components of preimages of $\bf F$. If the original Hamiltonian system is
nonresonant, then $\mathcal L$ is essentially unique, i.e. it does not depend on the
choice of the moment map $\bf F$.

Since we are dealing with topological aspects of integrable systems, we will be more
interested in the singular fibration ${\mathcal L}$ than in the moment map $\bf F$.
In particular, throughout this paper, we will adopt the following definition of
integrable systems:

\begin{define}
A singular Lagrangian torus fibration $ \pi : (M^{2n}, \omega, {\mathcal L}) \to O $
is called an {\it integrable system} (from the geometric point of view) if near each
fiber of $\mathcal L$ (i.e. preimage of $\pi$) it is defined by a set of $n$
commuting (with respect to the Poisson bracket) functionally independent functions
on $M^{2n}$. Two integrable systems are called {\it topologically equivalent} if
there is a fibration-preserving homeomorphism between them, and they are called {\it
symplectically equivalent} if there is a smooth fibration-preserving
symplectomorphism between them.
\end{define}

Notice that, in the above definition, we don't require the global existence of a
moment map. We just require it to exist in a neighborhood of every singular fiber.

The affine monodromy sheaf is a sheaf over the base space $O$ and is defined as
follows (see Subsection \ref{subsection:affinemonodromy}):

\begin{define}
The sheaf $\mathcal R$ over the base space $O$, which associates to each open subset
$U \subset O$ the free Abelian group $R(U)$ of symplectic system-preserving
${\mathbb S}^1$-actions in $(\pi^{-1}(U), \omega, {\mathcal L})$ is called the {\it
affine monodromy sheaf} of the system.
\end{define}

Besides affine monodromy, we will need another notion of monodromy which we call
{\it homological monodromy}, which involves first homology groups of the strata of
the fibers of the system and of the fibers themselves, and which will be explained
in Subsection \ref{subsection:rough}. In Subsection \ref{subsection:rough} we
introduce the notion of {\it rough equivalence} of integrable systems, which may be
reformulated as follows:

\begin{define}
Two integrable systems over the same base spaces are called {\it roughly
topologically equivalent} if they have the same singularities topologically, and the
same homological monodromy.
\end{define}

It is evident that if two systems are topologically equivalent, then they are also
roughly topologically equivalent, after an appropriate identification of their base
spaces. And under some hypotheses made in Section \ref{section:automorphisms} about
``genericity'' of systems under consideration, if two systems are roughly
topologically equivalent then they will have the same affine monodromy.

In order to define the Chern class of an integrable system, we will have to compare
it to a {\it reference system} which is roughly topologically equivalent to it. (For
regular fibrations, there is a natural choice of the reference system, which is the
one with a global section, so one does not have to mention it. But in our general
case there is no such a-priori choice). The definition of the Chern class involves
some cohomological exact sequence and is explained in Subsection
\ref{subsection:definitions}. The Chern class is an element of $H^2(O,{\mathcal
R})$, where $O$ is the base space and $\mathcal R$ is the affine monodromy sheaf.
Our main result is the following (see Subsection \ref{subsection:classification})

\begin{thm}
\label{thm:topological} Two roughly topologically equivalent integrable Hamiltonian
systems are topologically equivalent if and only if they have the same Chern class
with respect to a common reference system.
\end{thm}

What makes the above theorem effective is that in many cases it is relatively easy
to compute the cohomology group $H^2(O,{\mathcal R})$. Though ${\mathcal R}$ is not
a locally constant sheaf in general, it is locally constant on the strata of $O$,
and is a ``constructible'' free Abelian sheaf nevertheless.

Similarly, we have the following symplectic classification of integrable systems: in
Subsection \ref{subsection:rough} we introduce the notion of {\it rough symplectic
equivalence}, which is rough topological equivalence plus a condition of the
symplectic nature of the involved local automorphisms. In particular, if two systems
are roughly symplectically equivalent then their singularities are symplectically
equivalent (see e.g. \cite{DMT,Vu3,BF} for some results on symplectic invariants of
singularities of integrable systems). Then we introduce in Subsection
\ref{subsection:definitions} the {\it Lagrangian class}, which is a characteristic
class which lies in $H^1(O, {\mathcal Z}^{1}/{\mathcal R})$. Here ${\mathcal Z}^{1}$
is the sheaf of local closed differential 1-forms on $O$ (see Subsection
\ref{subsection:DeRham}), and there is a natural injection from $\mathcal R$ to
${\mathcal Z}^1$. We have (see Subsection \ref{subsection:classification}):

\begin{thm}
\label{thm:symplectic} Two roughly symplectically equivalent integrable Hamiltonian
systems are symplectically equivalent if and only if they have the same Lagrangian
class with respect to a common reference system.
\end{thm}

Theorem \ref{thm:symplectic} is in fact much easier to prove than Theorem
\ref{thm:topological}, because the sheaf of local system-preserving
symplectomorphisms is much smaller than the sheaf of local system-preserving
homeomorphisms. On the other hand, it is a highly non-trivial problem to classify
higher corank singularities of integrable systems symplectically. For systems
without singularities or with only elliptic singularities, Theorem
\ref{thm:symplectic} coincides with some results obtained earlier by Dazord and
Delzant \cite{DD} and Boucetta and Molino \cite{BM}.

A bonus of our classification results is the possibility to construct new integrable
systems and underlying symplectic structures from the old ones by means of surgery.
We call it {\it integrable surgery}, and give a few examples of this method (exotic
symplectic spaces, toric manifolds, K3, etc.) in Subsection
\ref{subsection:surgery}.

In this introduction, we often use the word ``generic'' without explaining what it
is. The problem is, though we know intuitively what does a ``generic'' integrable
system mean, we don't have a clear-cut definition, and we know very little about
degenerate higher corank singularities. So in Section \ref{section:automorphisms},
we will give a series of 7 ``very reasonable'' hypotheses, (H1)-(H7), about
singularities of integrable systems, and conjecture that all singularities of a
``generic'' integrable system, whatever it means, must satisfy these hypotheses.
These hypotheses have been verified for nondegenerate singularities. Theorem
\ref{thm:topological} is proved under the assumption that hypotheses (H1)-(H5) are
satisfied, and Theorem \ref{thm:symplectic} is proved under all 7 hypotheses.

The rest of this paper is organized as follows: In Section \ref{section:regular} we
give an exposition of the theory of regular Lagrangian torus fibrations developed by
Duistermaat, Dazord and Delzant \cite{Duistermaat,DD} (see also \cite{Marle,BM}).
Though most of the material of this section is not new, we include it here for the
convenience of readers, and to make it easier to see the similarities between the
regular case and the case with singularities. Section \ref{section:automorphisms}
contains the main preparation work of this paper, where we will study the structure
of the base space, local automorphisms of singularities, and write down a series of
hypotheses about singularities. In particular, we will define and study the
necessary sheaves for our characteristic classes there. Section
\ref{section:classification} contains the definition of characteristic classes, the
classification theorems, and a discussion of the realization problem and integrable
surgery.

Remark: In some of our previous papers, we used the words ``associated singular
Lagrangian {\it foliation}'' to call the associated fibration $\mathcal L$ of an
integrable system, but now we feel that the word {\it fibration} is a more correct
one. Likewise, we will use the word {\it fiber} instead of {\it leaf} when we talk
about a connected component of a level set of the moment map. This new convention
will save us from confusions when dealing with hyperbolic-type singular fibers, i.e.
fibers that contain more than one orbit of the Poisson action of the moment map,
because each (singular) orbit of the Poisson action is a (singular) leaf of the
associated singular foliation in the sense of Stefan-Sussmann.
\\

\section{Regular Lagrangian torus fibrations}
\label{section:regular}

If we throw out all singular fibers from an integrable system, then what remains is
a regular Lagrangian torus fibration.

\subsection{Action-angle coordinates}
\label{subsection:AA}

Let $\pi : (M^{2n}, \omega, {\mathcal L}) \to O$ be a regular Lagrangian fibration
with compact fibers. Then according to Arnold-Liouville theorem, each fiber of this
fibration is a Lagrangian torus of the symplectic manifold $(M^{2n}, \omega)$,
called a {\it Liouville torus}. Moreover, for each point $x \in O$ there is a
neighborhood $D^n = D(x)$ of $x$ in $O$ such that $(\pi^{-1}(D^n), \omega) \to D^n$
can be written as $(D^n \times {\Bbb T}^n, \sum_1^n dp_i \wedge dq_i) \to D^n$ via a
fibration-preserving symplectomorphism, where $(p_i)$ is a system of coordinates on
$D^n$ and ($q_i$ mod 1) is a system of periodic coordinates on ${\Bbb T}^n$. The
functions $p_i$ and $q_i$ are called {\it action} and {\it angle} coordinates,
respectively.

If $(u_i, v_i)$ is another system of action-angle
coordinates in $(\pi^{-1}(D^n), \omega, {\mathcal L})$, then we have
$$
\begin{pmatrix} u_1 \\ \vdots \\ u_n \end{pmatrix} = A
\begin{pmatrix} p_1 \\ \vdots \\ p_n \end{pmatrix} +
\begin{pmatrix} c_1 \\ \vdots \\ c_n \end{pmatrix} \; , \;
\begin{pmatrix} v_1 \\ \vdots \\ v_n \end{pmatrix} = (A^{-1})^T
\begin{pmatrix} q_1 \\ \vdots \\ q_n \end{pmatrix} +
\begin{pmatrix} g_1(p_i) \\ \vdots \\ g_n(p_i) \end{pmatrix},
$$
where $A$ is an element of $GL(n, {\Bbb Z})$, $c_i$ are
constants, and $\sum_1^n g_i dp_i$ is a closed differential 1-form on $D^n$.

In particular, the two local systems of action coordinates $(p_i)$ and $(u_i)$ on
$O$ are related by an integral affine transformation. Thus the base space $O$ admits
a unique natural integral affine structure. This integral affine structure provides
$O$ with a volume element, which is equal to $dp_1 \dots dp_n$ in any local system
of action coordinates, and the volume of $O$ is equal to the volume of $(M^{2n},
\omega)$ (for the standard volume form $\omega^n / n!$). Each Liouville torus also
admits a unique natural affine structure of a flat torus ${\Bbb T}^n = {\Bbb R}^n /
{\Bbb Z}^n$, given by a system of angle coordinates.

Given a regular Lagrangian torus fibration $\pi: (M^{2n}, \omega, {\mathcal L}) \to
O$, we can ask if there are global action-angle coordinates. That is, can $(M^{2n},
\omega) \to O$ be written in the form
$$
(O \times {\Bbb T}^n, \sum_1^n dp_i \wedge dq_i) \to O
$$
where $(p_i): O \to {\Bbb R}^n$ is an immersion and ($q_i$ mod 1) is a system of
periodic coordinates on ${\Bbb T}^n$.

A natural way to solve the above problem is via {\it obstruction theory}. If $\pi:
(M^{2n}, \omega) \to O$ admits a global system of action-angle coordinates, then it
has the following properties:

a) $\pi: M^{2n} \to O$ is a principal ${\Bbb T}^n$-bundle.

b) $\pi: M^{2n} \to O$ has a global section.

c) Moreover, it has a global {\it Lagrangian} section.

Conversely, if the above conditions are satisfied then one can show easily that
$\pi: (M^{2n}, \omega) \to O$ admits global action-angle coordinates.

The obstruction for the condition a) to be fulfilled will be called the {\it
monodromy}, or also the {\it affine monodromy}, because it can be determined
completely by the affine structure of the base space $O$. Obstructions to b) and c)
will be characterized by the so-called {\it Chern class} and {\it Lagrangian class},
respectively.

\subsection{Monodromy}
\label{subsection:monodromy}

Given a Lagrangian torus fibration $\pi: (M^{2n}, \omega) \to O$, we will associate
to it the ${\Bbb Z}^n$-bundle of first homology groups of the fibers of ${\mathcal
L}$, denoted by $ E_{\Bbb Z} \stackrel{H_1({\Bbb T}^n, {\Bbb Z})}{\longrightarrow}
O.$ The holonomy of this bundle and is an element of $\hom(\pi_1(O), GL(n, {\Bbb
Z}))$, defined up to conjugacy, and is called the {\it monodromy} of the torus
fibration.

The symplectic form $\omega$ gives rise to a natural isomorphism from the vector
bundle $ E_{\Bbb R} \stackrel{H_1({\Bbb T}^n, {\Bbb R})}{\longrightarrow} O$ (of
first homology groups over $\Bbb R$ of the fibers of $\mathcal L$) to the cotangent
bundle $T^{\ast}O$, defined as follows: Let $T_x = \pi^{-1}(x)$ be a fiber of
$\mathcal L$. Then $T_x$ has a unique canonical flat structure, and constant vector
fields on $T_x$ can be identified with $H_1(T_x, {\Bbb R})$ via ``rotation numbers".
On the other hand, if $X$ is a constant vector field on $T_x$, then the covector
field $\alpha(X) = - i_X \omega$ is the pull-back of a covector $\alpha$ on $O$ at
$x$, i.e. an element of $T^{\ast}O$.

Notice that $ E_{\Bbb Z} \stackrel{H_1({\Bbb T}^n, {\Bbb Z})}{\longrightarrow} O $
is a discrete subbundle of $ E_{\Bbb R} \stackrel{H_1({\Bbb T}^n, {\Bbb
R})}{\longrightarrow} O $. Under the aforementioned identification of $E_{\Bbb R}$
with $T^{\ast}O$, $E_{\Bbb Z}$ maps to a discrete subbundle of $T^{\ast}O$,
consisting of ``integral'' covectors. We will denoted this subbundle, or the
discrete sheaf associated to it, by $\mathcal R$. It follows from Arnold-Liouville
theorem that local sections of $\mathcal R$ are local differential 1-forms on $O$
which can be written as $\sum m_i dp_i$ in some local system of action coordinates
$(p_i)$, with $m_i \in {\Bbb Z}$. Thus $\mathcal R$ can be completely determined by
the integral affine structure of $O$. We will call $\mathcal R$ the {\it affine
monodromy sheaf}. Since $E_{\Bbb Z}$ is isomorphic to $\mathcal R$, the monodromy of
the system is completely determined by the integral affine structure of $O$, and
will also be called the {\it affine monodromy}.

There is another characterization of $\mathcal R$ as follows: Arnold-Liouville
theorem implies that each local differential 1-form on $O$ of the type $\sum m_i
dp_i$ with $m_i \in {\Bbb Z}$ in a local system of action coordinates $(p_i)$ gives
rise to a symplectic vector field $\sum m_i X_{p_i}$ which generates a symplectic
${\Bbb S}^1$-action which preserves the system, and vice versa. Thus $\mathcal R$ is
isomorphic to the sheaf of local system-preserving symplectic ${\Bbb S}^1$-actions.

First examples of integrable systems with nontrivial monodromy, namely
the spherical pendulum and the Lagrange top, were observed by Cushman and
others (e.g., \cite{CK,Duistermaat}). In these examples and all other
known examples arising from classical mechanics and physics,
the nontriviality of the monodromy is due to the
presence of the so-called focus-focus singularities
(see e.g. \cite{ZungFocus}).

\subsection{Chern and Lagrangian classes}

The Chern class can be defined as the obstruction for the torus
fibration $M \to O$ to admit a global section.
Let $(U_i)$ be a trivializing open covering of $O$. Over each $U_i$ there
is a smooth section, denoted by $s_i$. The difference between two local
sections $s_i$ and $s_j$, over $U_i \cap U_j$, can be written as
$$
\mu_{ij} = s_j - s_i \in C^{\infty}(E_{\Bbb R}/ E_{\Bbb Z}) (U_i \cap U_j) \cong
C^{\infty}(T^{\ast}O / {\mathcal R}) (U_i \cap U_j) .
$$
Here $C^{\infty}(.)$ denotes the sheaf of smooth sections, and $E_{\Bbb R}$ and
$E_{\Bbb Z}$ are the first cohomology bundles defined in the previous subsection. It
is immediate that $(\mu_{ij})$ is a 1-cocycle, and it defines a \v{C}ech first
cohomology class, not depending on the choice of sections:
$$
\hat{\mu} \in H^1(O, C^{\infty}(T^{\ast}O / {\mathcal R})) .
$$
Since $C^{\infty}(T^{\ast}O)$ is a fine sheaf, from the short exact sequence
$$
0 \to {\mathcal R} \to C^{\infty}(T^{\ast}O) \to
C^{\infty}(T^{\ast}O/{\mathcal R})) \to 0
$$
we obtain that the coboundary map $\delta :
H^1(O, C^{\infty}(T^{\ast}O/{\mathcal R})) \to H^2(O, {\mathcal R})$
in the associated long exact sequence is an isomorphism.

The image $\mu_{C}$ of $\hat{\mu}$ in $H^2(O, {\mathcal R})$ under the isomorphism
$\delta$ is called the {\it Chern class} \cite{Duistermaat}. In the case of trivial
monodromy,  $\mu_{C}$ coincides with the usual Chern class of principal ${\Bbb T}^n$
bundles (cf. \cite{DD}).

If one requires local sections $s_i$ to be Lagrangian, then one has that
$$
\mu_{ij} \in {\mathcal Z}(T^{\ast}O / {\mathcal R}) (U_i \cap U_j)
$$
where ${\mathcal Z}$ means closed 1-forms, and it will define another
cohomology class which we will call the {\it Lagrangian class}:
$$
\mu_{L} \in H^1(O, {\mathcal Z}(T^{\ast}O / {\mathcal R}))
$$

There is another short exact sequence
$$
0 \to {\mathcal R} \to {\mathcal Z}(T^{\ast} O)
\to {\mathcal Z}(T^{\ast} O /{\mathcal R}) \to 0 ,
$$
which leads to the following long exact sequence
$$\begin{array}{c}
\ldots \to H^1(O, {\mathcal R}) \stackrel{\hat{d}}{\to}
H^1(O, {\mathcal Z}(T^{\ast}O)) \equiv H^2(O, {\Bbb R})
\to H^1(O, {\mathcal Z}(T^{\ast}O /{\mathcal R})) \\
\stackrel{\Delta}{\to} H^2(O, {\mathcal R})
\stackrel{\hat{d}}{\to} H^2(O, {\mathcal Z}(T^{\ast}O_0))
\equiv H^3(O, {\Bbb R}) \to H^2(O, {\mathcal Z}(T^{\ast}O/{\mathcal R})) \to \ldots
\end{array}
$$

Under the maps $\Delta$ and $\hat{d}$ we have
$\mu_{L} \stackrel{\Delta}{\mapsto} \mu_{C} \stackrel{\hat{d}}{\mapsto} 0.$

Thus, if the integral affine manifold $O$ is given, then any element of
the first cohomology group
$H^1(O, {\mathcal Z}(T^{\ast}O/{\mathcal R}))$ will be the Lagrangian
class of some Lagrangian torus fibration over $O$, and the necessary and
sufficient condition for an element $\mu$ in $H^2(O, {\mathcal R})$ to be the
Chern class of some Lagrangian torus fibration is that
$\hat{d}(\mu) = 0$. To each element $\mu_{C} \in H^2(O, {\mathcal R})$
such that $\hat{d}(\mu_{C}) = 0$, there are
$H^2(O, {\Bbb R}) / \hat{d}H^1(O, {\mathcal R})$ choices of the element
$\mu_{L}$ such that $\Delta(\mu_{L}) = \mu_{C}$, and each choice
corresponds to a symplectically different Lagrangian torus fibration
with the same Chern class $\mu_C$. If $\mu_{L} = 0$ then the
corresponding fibration is symplectically equivalent to
$T^{\ast}O / {\mathcal R} \longrightarrow O$ (cf. \cite{DD}).

In particular, if the base space $O$ is 2-connected :
$\pi_1(O) = \pi_2 (O) = 0$, then there is no room for the monodromy
and the Lagrangian class, so there always exists a global system of
action-angle coordinates, as first observed by
Nekhoroshev \cite{Nekhoroshev}.

It is clear from the definition of the Chern class that if two regular Lagrangian
torus fibrations $(M_1, \omega_1, {\mathcal L}_1) \to O$ and $(M_2, \omega_2,
{\mathcal L}_2) \to O$ over the same base space $O$ admit a diffeomorphism $\phi :
M_1 \to M_2$ which projects to the identity map on $O$ and preserves the flat
structure of each torus, then they have the same Chern class. In fact, the Chern
class is a topological invariant, in the sense that if $\phi : M_1 \to M_2$ is a
homeomorphism which projects to the identity map on $O$, but which need not be a
diffeomorphism and need not preserve the flat structure of the fibers, then the two
fibrations still have the same Chern class. In the following sections we will show
this fact for the more general case of systems with singularities.

\begin{example}{Systems over flat tori and Kodaira-Thurston example}
\end{example}

Assume that $O = {\Bbb T}^2 = {\Bbb R}^2 / {\Bbb Z}^2$, the standard flat torus with
the integral affine structure induced from ${\Bbb R}^2$. Then $H^3(O, {\Bbb R}) =
0$, and every element $\mu_{DC} \in H^2(O, {\mathcal R}) = {\Bbb Z}^2$ is realizable
as the Chern class of some Lagrangian torus fibration over $O$. The automorphism
group of the base space acts on $H^2(O, {\mathcal R})$, and the quotient space is
isomorphic to ${\Bbb Z}_+$ (nonnegative integers). Thus each integrable system with
the base space ${\Bbb T}^2$ is characterized topologically by a nonnegative integer
$m$, and its ambient symplectic manifold $M^4_m$ has $H_1(M^4_m, {\Bbb Z}) = {\Bbb
Z}^3 \oplus ({\Bbb Z} / m {\Bbb Z})$ as can be computed easily. For each $m$ there
are $H^2(O, {\Bbb R}) / \hat{d} H^1 (O, {\mathcal R}) = {\Bbb R} / {\Bbb Z}$ choices
of the symplectic structure on the fibration $M^4_m \to O$, up to symplectic
equivalence. Notice that the fibrations $M^4_m \to {\Bbb T}^2$ are topologically the
same as a series of elliptic fibrations over an elliptic curve (see e.g.
\cite{BPV}). In particular, when $m=1$, $M^4_1$ is the well-known Kodaira-Thurston
example (see e.g. \cite{MS}) of a manifold admitting both a complex and a symplectic
structure but not a K\"ahler structure.

If, for example, $O = {\Bbb R}^k / \Sigma$ where $\Sigma$ is a lattice of ${\Bbb
R}^k$, and $k \geq 3$, then not every element $\mu$ of $H^2(O, {\mathcal R})$ will
satisfy the condition $\hat{d}(\mu) = 0$. If the lattice $\Sigma$ is irrational,
then it may happen that the operator $\hat{d}$ is injective, and the only Lagrangian
torus fibration over $O$ is the one which admits a Lagrangian section. \\


\section{Base space and sheaves of local automorphisms}
\label{section:automorphisms}

From now on, $\pi : (M, \omega, {\mathcal L}) \to O$ will always denote an
integrable system which may admit singularities.

\subsection{Singularities of integrable systems}

By a singularity of $\pi : (M, \omega, {\mathcal L}) \to O$ we will mean the germ of
the fibration ${\mathcal L}$ at a singular fiber $N_x = \pi^{-1}(x)$, $x \in O$, and
will denote it by $\pi : (U(N_x), \omega, {\mathcal L}) \to (U(x))$, where $U(N_x) =
\pi^{-1}(U(x))$ is a tubular neighborhood of $N_x$. Two singularities are called
{\it topologically} (resp. {\it symplectically}) {\it equivalent} if their fibration
germs are homeomorphic (resp. symplectomorphic).

The {\it rank} of a singular fiber $N$ of $\mathcal L$ is ${\rm rank\ } N = \max_{p
\in N} {\rm rank\ } p$, where, by definition,
$$
{\rm rank\ } p = \max_{\bf F} \dim <dF_1(p),\dots,dF_n(p)>,
$$
where the maximum is taken over all possible moment maps
${\bf F} = (F_1,...,F_n) : (M,\omega, {\mathcal L}) \to {\Bbb R}^n$,
and $<dF_1(p),\dots,dF_n(p)>$ denotes the linear span
of the covectors $dF_1(p)$,...,$dF_n(p)$ in $T^{\ast}_p M$.

We put ${\rm corank\ } N = n - {\rm rank}\ N$ and ${\rm corank}\ p = n - {\rm rank\ } p$, where
$n = 1/2 \dim M$. If ${\rm rank\ } p < n$, then $p$ is called a {\it singular
point} of the system. If  ${\rm rank\ } p = 0$ then $p$ is called a {\it fixed
point}. The rank and corank of a singularity
$\pi : (U(N_x), \omega, {\mathcal L}) \to U(x)$ is, by definition,
the rank and corank of $N_x$.

A fixed point $p$ of the system is called {\it nondegenerate} if it satisfies the
following condition: there is a moment map $(F_1,\dots,F_n): (M, p, \omega,
{\mathcal L}) \to ({\Bbb R}^n, 0)$ such that the quadratic parts
$F^{(2)}_1,\dots,F^{(2)}_n$ of $F_1\dots,F_n$ in a symplectic system of coordinates
at $p$ will form a Cartan subalgebra of the Lie algebra of quadratic functions on
${\Bbb R}^{2n}$ under the standard Poisson bracket. Recall that the algebra of
quadratic functions on $({\Bbb R}^{2n}, \omega_0)$ is naturally isomorphic to the
simple Lie algebra $sp(2n, {\Bbb R})$, the functions $F^{(2)}_1,\dots,F^{(2)}_n$
Poisson-commute (because $F_1,\dots,F_n$ Poisson-commute) and span an Abelian
subalgebra of $sp(2n, {\Bbb R})$. The nondegeneracy condition means that this
subalgebra is of dimension $n$ and consists of semi-simple elements, i.e. it is a
Cartan subalgebra. A classical theorem of Williamson \cite{Williamson} (which
essentially classifies Cartan subalgebras of $sp(2n, {\Bbb R})$ up to conjugacy)
implies that, for a nondegenerate fixed point $p$, there is a moment map
$(F_1,\dots,F_n)$ whose quadratic part at $p$ can be decomposed into components of 3
types : elliptic ($F^{(2)}_i = p_i^2 + q_i^2$), hyperbolic ($F^{(2)}_i = p_iq_i$),
and focus-focus ($F^{(2)}_i = p_iq_i + p_{i+1}q_{i+1}, F^{(2)}_{i+1} = p_iq_{i+1} -
p_{i+1}q_i$), in a symplectic system of coordinates $(p_i,q_i)$ (the symplectic form
is $\omega = \sum dp_i \wedge dq_i$). Note that each focus-focus component consists
of two functions. If there are $k_e$ elliptic, $k_h$ hyperbolic and $k_f$
focus-focus components $(k_e + k_h + 2 k_f = n)$, then we will say that the {\it
Williamson type} of $p$ is $(k_e,k_h,k_f)$ (cf. \cite{ZungAL}). The local normal
form theorem for nondegenerate fixed points of integrable systems says that an
integrable Hamiltonian system near a nondegenerate fixed point is locally
symplectically equivalent to a system given by a quadratic moment map on the
standard symplectic vector space $({\Bbb R}^{2n}, \omega_0)$. This normal form
theorem has been proved in the analytic case by R\"ussmann \cite{Russmann} and Vey
\cite{Vey} (see also \cite{ZungNormal}).  In the smooth case, it has been proved
partially by Colin de Verdier and Vey \cite{CV}, Eliasson \cite{Eliasson}, and
Dufour and Molino \cite{DM}. (We don't know yet of any complete proof for the smooth
non-elliptic case, except for the case with only one degree of freedom).

A singular point $p$ with ${\rm rank\ } p = r > 0$ is called {\it nondegenerate} if
it becomes a nondegenerate fixed point for a local integrable Hamiltonian system
with $n-r$ degrees of freedom after a local Marsden-Weinstein reduction. A
singularity $\pi : (U(N_x), \omega, {\mathcal L}) \to U(x)$ is called {\it
nondegenerate} if all singular points in $N_x$ are nondegenerate, plus a natural
additional condition which is called ``topological stability'' in \cite{ZungAL}.

Remark: In some references the above additional condition is not included in the
definition of nondegenerate singularities, but we will include in here so that we
can use the decomposition theorem mentioned below. In \cite{BF}, Bolsinov and
Fomenko suggested another name for this additional condition, which looks better and
more to the point than our original name ``topological stability'', but
unfortunately I could not translate this name into English.

In \cite{ZungAL} we studied nondegenerate singularities of integrable Hamiltonian
systems, where we showed, among other things, that they are topologically equivalent
to almost direct products of simplest singularities. More precisely, if $(U(N_x),
{\mathcal L})$ is a nondegenerate singularity, then it is homeomorphic to
$$
\begin{array}{c}
\{ (U(T^r), {\mathcal L_r}) \times (P^2(N_1), {\mathcal L}_1) \times \dots
\times (P^2(N_{k_e + k_h}), {\mathcal L}_{k_e + k_h}) \times \\
\times (P^4(N_1'), {\mathcal L}_1') \times \dots
\times (P^4(N_{k_f}'), {\mathcal L}_{k_f}') \}/ {\Gamma}
\end{array}
$$
where $(U(T^r), {\mathcal L_r})$ denotes a regular system with $r$ degrees of
freedom near a regular torus, $(P^2(N_i), {\mathcal L}_i)$ for $1 \leq i \leq k_e +
k_h$ denotes a corank-1 nondegenerate elliptic or hyperbolic singularity of an
system with 1 degree of freedom, $(P^4(N_i'), {\mathcal L}_i')$ for $1 \leq i \leq
k_f$ denotes a focus-focus singularity of a system with 2 degrees of freedom,
$\Gamma$ is a finite group that acts on the above product freely and component-wise.
Moreover, it acts trivially on all possible elliptic components of the product
(there are $k_e$ such components if $k_e > 0$). Nondegenerate corank-1 singularities
are called ``atoms'' in the works of Fomenko and his collaborators (see e.g.
\cite{BF,Fomenko5}). Focus-focus corank-2 singularities are classified topologically
in \cite{ZungFocus}.

To our knowledge, most singularities of integrable Hamiltonian systems
are nondegenerate. But starting with 2 degrees of freedom, there are also
degenerate singularities. The situation is similar to that of smooth
functions: most singularities of smooth functions are nondegenerate,
but there are degenerate singularities whose miniversal deformation is
of finite dimension $k$ and which can appear in a generic way in
$k$-dimensional families of functions (see e.g. \cite{ZungCorank1}).

We will make a series of hypotheses about singularities of ``generic'' integrable
Hamiltonian systems. All nondegenerate singularities will satisfy these hypotheses.
We believe that ``generic'' degenerate singularities will also satisfy these
hypotheses. The first hypothesis is:

\begin{quote}
{\bf (H1)}{\it Each fiber $N$ of ${\mathcal L}$ is a disjoint union of a finite
number of submanifolds $N_i$ such that if $p \in N_i$ then ${\rm rank\ } p = \dim
N_i$.}
\end{quote}

Recall that if $p \in N_x \subset M$ has rank $r$, then there is a moment map ${\bf
F} = (F_1,\dots,F_n): (M,\omega) \to {\Bbb R}^n$ such that the rank of ${\bf F}$ at
$p$ is $r$, and the orbit of the Poisson ${\Bbb R}^n$-action generated by ${\bf F}$
which contains $p$ is an immersed $r$-dimensional submanifold contained in $N_x$.
Hypothesis (H1) says that each such orbit is a stratum of $N_x$ in a natural sense.
In particular, $\dim N_x \leq n$ for any fiber $N_x$ of $\mathcal L$. If, for
example,$(U(N_x), \omega, {\mathcal L})$ is a nondegenerate singularity of
Williamson type $(k_e,k_h,k_f)$ with $k_e > 0$ (where $k_e$ is the number of
elliptic components), then $\dim N_x = n - k_e < n$.

\subsection{Stratification of the base space}
\label{subsection:stratification}

The base space $O$ of an integrable system $\pi : (M, \omega, {\mathcal L}) \to O$
has a natural topology induced from $M$. We will always assume $O$ to be separated
and paracompact (it follows from the existence of a proper moment map). If $\mathcal
L$ is regular then $O$ is a manifold. If $\mathcal L$ has only nondegenerate
elliptic singularities then $O$ is a manifold with corners (see e.g. \cite{BM}). In
the general case, $O$ is not a manifold, but we can try to give it a natural
stratification.

In this paper, a separated paracompact topological space $Q$ will be called
a {\it stratified manifold} of dimension $n$ if $Q_i$ can be written as a
disjoint union of topological manifolds $Q_i$, called strata of $Q$, in
such a way that:

a) $\max \dim Q_i = n$

b) For each $i$, the boundary of the stratum $Q_i$ is a union of strata
of dimension smaller than $\dim Q_i$.

c) If $\dim Q_i = k$ then for each $x \in Q_i$ there is a neighborhood
$U(x)$ of $x$ in $V$ which is homeomorphic to the direct product of a
$k$-dimensional disk $D^k$ with a cone over a stratified
$(n-k-1)$-dimensional manifold with a finite number of strata. Such a
neighborhood will be called a {\it standard}, or {\it star-shaped},
neighborhood of $x$.

Given a singular point $x$ of the base space $O$, we will denote by $S_x$ the
connected component of the set of all points $y \in O$ such that the singularity at
$y$ (i.e. at the singular fiber $N_y = \pi^{-1}(y)$ of the system) is topologically
equivalent to the singularity at $x$. If $x$ is regular then $S_x$ is a connected
component of the regular part of $O$. Our next hypothesis is:

\begin{quote}
{\bf (H2)} {\it The base space $O$ is a stratified manifold for which
each $S_x$ defined above is a stratum. If $U(x) \in O$ is a star-shaped
neighborhood of a point $x \in O$, and $\phi : U(x) \to D^k \times V(x)$
is a homeomorphism, where $D^k$ is a $k$-dimensional disk and
$V(x) \in O$ is a local stratified submanifold transversal to $S_x$ at
$x$ (so $V(x)$ is homeomorphic to a cone over a $(n-k-1)$-dimensional
stratified manifold), then there is a homeomorphism
$\Phi : \pi^{-1}(U(x)) \to D^k \times \pi^{-1}(V(x))$ which
makes the following diagram commutative:
$$
\begin{array}{ccc}
\pi^{-1}(U(x)) & \stackrel{\Phi}{\longrightarrow} & D^k \times \pi^{-1}(V(x)) \\
\pi \downarrow & & \downarrow (id, \pi) \\
U(x) & \stackrel{\phi}{\longrightarrow} & D^k \times V(x)
\end{array}
$$
}
\end{quote}

It is clear that Hypothesis (H2) has a local character : it is satisfied if it is
satisfied locally near every singular point. Hypothesis (H2) justifies the use of a
tubular neighborhood of a singular fiber to denote a singularity. In fact, it
follows from this hypothesis that if $U_1(x)$ and $U_2(x)$ are two small star-shaped
neighborhood of a point $x$ in $O$, then the singular fibration
$(\pi^{-1}(U_1(x)),{\mathcal L})$ is topologically equivalent to the singular
fibration $(\pi^{-1}(U_2(x)),{\mathcal L})$.

Let $O = \sqcup S_x$ be the stratification of the base space $O$ as above. Then
we can  replace $S_i$ by {\it thickened strata} as follows:
Order the strata of $O$ by dimension: $S_1, S_2, S_3, ...$,
with $\dim S_1 \leq \dim S_2 \leq \dim S_3 \leq ...$. Denote by $\bar{S_1}$ a
small closed tubular neighborhood of $S_1$ in $O$,  $\bar{S_2}$ a sufficiently
small closed tubular neighborhood of $(O \setminus \bar{S_1}) \cap S_2$ in
$(O \setminus \bar{S_1})$, $\bar{S_3}$ a sufficiently
small closed tubular neighborhood of $(O \setminus (\bar{S_1} \cup \bar{S_2})) \cap S_3$ in
$(O \setminus (\bar{S_1} \cup \bar{S_2}))$, etc. Then we have $O = \sqcup \bar{S_i}$.
This decomposition of $O$ into the disjoint union of $\bar{S_i}$ is unique up
to homeomorphisms, and is called a {\it thickened stratification} of $O$.
The sets $\bar{S_i}$ are called {\it thickened strata} of $O$.

Similarly, if $U(x)$ is a star-shaped neighborhood of a point $x \in O$, then
we also have a thickened stratification of $U(x)$, $U(x) = \sqcup \bar{U_i}$,
which consists of a finite number of thickened strata $\bar{U_i}$ and
is unique up to homeomorphisms. In an appropriate thickened stratification of
$O$, $O =\sqcup \bar{S_i}$, we can put
$\bar{U_i} = U(x) \cap \bar{S_i}$ (and then throw out those
$\bar{U_i} = U(x) \cap \bar{S_i}$ which are empty).

\subsection{Local ${\Bbb S}^1$-actions and affine monodromy}
\label{subsection:affinemonodromy}

Let $U$ be an open subset of the base space $O$. We will denote by $R(U)$ the set of
symplectic ${\Bbb S}^1$-actions on $(\pi^{-1}(U), \omega)$ which preserve the system
(i.e. preserve every fiber of $\mathcal L$). We have the following obvious lemma:

\begin{lem}
$R(U)$ is a free Abelian group of rank less or equal to $n$. If $U_1$ is an
open subset of $U$ then there is a natural injection from $R(U)$ to $R(U_1)$.
If $\alpha$ is a system-preserving symplectic ${\Bbb S}^1$-action on $\pi^{-1}(U_1)$
and $\beta$ is a system-preserving symplectic ${\Bbb S}^1$-action on $\pi^{-1}(U_2)$
such that their restrictions to $\pi^{-1}(U_1 \bigcap U_2)$ are the same, then
there is a system-preserving symplectic ${\Bbb S}^1$-action on $\pi^{-1}(U_1 \bigcup U_2)$
which restricts to $\alpha$ and $\beta$ on $\pi^{-1}(U_1)$ and
$\pi^{-1}(U_2)$ respectively.
\end{lem}

Thus, the groups $R(U)$, $U \subset O$, form a free Abelian sheaf over the base
space $O$. We will denote this sheaf by $\mathcal R$, and call it the {\it affine
monodromy sheaf}, in analogy with the case of regular Lagrangian torus fibrations.

If $(g_1,\dots,g_m)$ is a basis of $R(U) \cong {\Bbb Z}^m$, then these ${\Bbb
S}^1$-actions $g_1,\dots,g_m$ commute, and together they generate a
system-preserving symplectic ${\Bbb T}^m$-action on $(\pi^{-1}(U), \omega, {\mathcal
L})$ which is free almost everywhere. Conversely if there is a system-preserving
symplectic action of a torus ${\Bbb T}^m$ on $(\pi^{-1}(U), \omega, {\mathcal L})$
which is locally free somewhere, then the composition of this ${\Bbb T}^m$-action
with homomorphisms from ${\Bbb S}^1$ to ${\Bbb T}^m$ gives rise to a subgroup of
$R(U)$ which is isomorphic to ${\Bbb Z}^m$. The classical Arnold-Liouville theorem
is essentially equivalent to the fact that, if $U$ is a disk in the regular part of
$O$, then $R(U)$ is isomorphic to ${\Bbb Z}^n$. In \cite{ZungAL} we have shown that
if $\pi : (\pi^{-1}(U(x)), \omega, {\mathcal L}) \to U(x)$ is a nondegenerate
singularity of rank $r$ and Williamson type $(k_e,k_h,k_f)$ ($r + k_e + k_h + 2 k_f
= n$), then $R(U(x))$ is isomorphic to ${\Bbb Z}^{r + k_e + k_f}$. In
\cite{ZungCorank1} we have shown that if $\pi : (\pi^{-1}(U(x)), \omega, {\mathcal
L}) \to U(x)$ is a degenerate singularity of corank $1$, then under some mild
conditions $R(U(x))$ will be isomorphic to ${\Bbb Z}^{n-1}$. Another simple result
which can be proved by the methods of \cite{ZungAL,ZungCorank1} is the following:

\begin{lem}
\label{lem:hyperbolic}
If $\pi : (\pi^{-1}(U(x)), \omega, {\mathcal L}) \to U(x)$ is a
singularity of rank $r$ such that
$\dim \pi^{-1}(x) = n = 1/2 \dim \pi^{-1}(U(x))$, then there is a locally-free
system-preserving symplectic action of ${\Bbb T}^r$ on
$(\pi^{-1}(U(x)), \omega, {\mathcal L})$. In particular,
the rank of $R(U(x))$ is greater or equal to $r$.
\end{lem}

{\it Sketch of the Proof}. Denote by $P$ an $n$-dimensional stratum of
$\pi^{-1}(x)$. It follows from the facts that $P$ is an orbit of the Poisson action
of a moment map $(F-1,\dots,F_n): (\pi^{-1}(U(x)), \omega, {\mathcal L}) \to {\Bbb
R}^n $ and ${\rm rank}\ \pi^{-1}(x) = r$, that $P$ is of the type ${\Bbb T}^k \times
{\Bbb R}^{n-k}$ with $k \geq r$. For each element $\gamma$ of the fundamental group
of $P$, there are numbers $a_1,\dots,a_n \in {\Bbb R}$ such that the vector field
$\sum a_i X_{F_i}$ is periodic of period 1 on $P$ and its orbits on $P$ are
homotopic to $\gamma$. There is a unique way to extend $a_i$ into smooth functions
on $\pi^{-1}(U(x))$ which are constant on each fiber of ${\mathcal L}$, such that
the vector field $\sum a_i X_{F_i}$ is periodic on $\pi^{-1}(U(x))$.
Arnold-Liouville theorem (for the regular part of $\pi^{-1}(U(x))$) assures that the
vector field $\sum a_i X_{F_i}$ is symplectic. It follows that there is a
system-preserving symplectic action of ${\Bbb T}^k$ on $(\pi^{-1}(U(x)), \omega,
{\mathcal L})$ which is free almost everywhere. We can choose a subgroup ${\Bbb T}^r
\subset {\Bbb T}^k$ such that the action of ${\Bbb T}^r$ will be locally-free
everywhere in $\pi^{-1}(U(x))$. $\blacksquare$

Let $ g : {\Bbb S}^1 \times \pi^{-1}(U) \longrightarrow \pi^{-1}(U) $ be a
system-preserving symplectic ${\Bbb S}^1$-action, where $U$ is an open subset of the
base space $O$. Then for each stratum $P$ of a fiber in $\pi^{-1}(U)$, $g$ preserves
$P$ and the orbits of $g$ on $P$ defines an element $\gamma_P$ in the fundamental
group of $P$. The association $P \mapsto \gamma_P \in \pi_1(P)$ is continuous in the
following sense: each continuous curve $\psi : [0,1] \to \pi^{-1}(U)$ can be
extended to a continuous family of loops $\phi: [0,1] \times S^1 \to \pi^{-1}(U)$ \
($\phi\mid_{[0,1] \times \{ e\}} = \psi$ where $e$ is a fixed element of $S^1$) such
that each loop $\phi\mid_{\{t\} \times S^1}$ lies on some stratum $P(t)$ of some
fiber of $\pi^{-1}(U)$ and is homotopic to $\gamma_{P(t)}$ on $P(t)$. It is evident
that the set of continuous association $\{ P \mapsto \gamma_P \in \pi_1(P) \}$ is an
Abelian group, which we will denote by $\Pi (U)$, and there is a unique natural
injective homomorphism from the group $R(U)$ to $\Pi(U)$. Our next hypothesis is:

\begin{quote}
{\bf (H3)} {\it If $\pi : (\pi^{-1}(U(x)), \omega, {\mathcal L}) \to U(x)$
is a singularity of rank $r$, then $R(U(x))$ is isomorphic to $\Pi(U(x))$,
and there is a locally-free system-preserving symplectic ${\Bbb T}^r$-action
on $(\pi^{-1}(U(x)), \omega, {\mathcal L})$. In particular,
${\rm rank}\ R(U(x)) = {\rm rank}\ \Pi(U(x)) \geq {\rm rank}\ \pi^{-1}(x)$.
}
\end{quote}

Since $\Pi(U(x))$ is clearly a topological invariant of
$(\pi^{-1}(U(x)), {\mathcal L})$, Hypothesis (H3) implies in particular that
the affine monodromy sheaf ${\mathcal R}$ is a topological invariant of the
system. Moreover, $R(U)$ is naturally isomorphic to $\Pi(U)$ for each
open subset $U$ of $O$.

\subsection{Affine structure of the base space}
\label{subsection:affine}

A local function $f$ on the base space $O$ is called an {\it integral affine
function} if its pull-back $f \circ \pi$ is a local smooth function on $M$ such that
the Hamiltonian vector field of $f \circ \pi$ is periodic of period $1/k$ for some
natural number $k$. In other words, local integral affine functions are functions
that generate local system-preserving ${\mathbb S}^1$-actions. Denote the sheaf of
local integral affine functions on $O$ by ${\mathcal I}$. Then we have the following
natural short exact sequence:
$$
0 \longrightarrow {\Bbb R} \longrightarrow {\mathcal I} \longrightarrow {\mathcal R}
\longrightarrow 0.
$$

We may view the sheaf ${\mathcal I}$ as the integral affine structure of $O$. A
local function on $O$ is called {\it affine} if it can be written as a linear
combination of integral affine functions.

As a side note, we remark that the fact that $O$ has a integral affine structure
imposes some conditions on the topology of $O$. For example, for 2-dimensional base
spaces we have the following analog of a Theorem of Milnor \cite{Milnor} concerning
affine structures on 2-surfaces:

Let $O$ be the base space of an integrable system with 2 degrees of freedom, whose
singularities are all nondegenerate. A subset $C$ of $O$ will be called a {\it
topological 2-stratum} of $O$ if $C$ is the union of a 2-dimensional stratum of $O$
(with respect to the stratification given in Subsection
\ref{subsection:stratification}) with all possible focus-focus singular points lying
on its boundary.

\begin{prop}
\label{prop:2domain} Let $C$ be a topological 2-stratum of the base $O$ of an
integrable Hamiltonian system on a compact 4-dimensional symplectic manifold $M^4$
(maybe with boundary). Assume that the system contains only nondegenerate
singularities, and the image of the boundary of $M^4$ under the projection to $O$
does not intersect with the closure of $C$ (if $M^4$ is closed then this condition
is satisfied automatically). Then $C$ is homeomorphic to either an annulus, a Mobius
band, a Klein bottle, a torus, a disk, a projective space, or a sphere (in case of
sphere or projective space, $C$ must contain focus-focus points).
\end{prop}

The proof of the above proposition is elementary and is left to the reader.

\subsection{Local homeomorphisms}

For each open subset $U \subset O$, we will denote by $A_t(U)$ the group of
homeomorphisms $\phi$ from $(\pi^{-1}(U), {\mathcal L})$ onto itself which satisfy
the following topological properties: $\phi$ preserves each stratum of each fiber of
$(\pi^{-1}(U), {\mathcal L})$ and induces the identity map on the fundamental group
of each stratum, and $\phi$ induces the identity map on the first homology group
(with integral coefficients) of each fiber of $(\pi^{-1}(U), {\mathcal L})$. The
association $U \mapsto A_t(U)$ is a non-Abelian sheaf over the base space $O$, which
will be denoted by ${\mathcal A}_t$ and called the {\it sheaf of admissible local
topological automorphisms} of the system. (Here {\it t} stands for ``topological'').

If $g : [0,1] \to A_t(U)$ is a continuous loop in $A_t(U)$, then its orbits on the
strata of the  fibers of $(\pi^{-1}(U), {\mathcal L})$ define an element of $\Pi(U)$
($ = R(U)$ by Hypothesis (H3)). Moreover, it follows from the definition of $A_t(U)$
that every element of $A_t(U)$ induces the identity map on $\Pi(U)$. Thus there is a
natural central extension $\hat{A}_t'(U)$ of $A_t(U)$ by $R(U)$: $0 \to R(U) \to
\hat{A}_t'(U) \to A_t(U) \to 0$. Topologically (with respect to the open-compact
topology), $ \hat{A}_t'(U)$ is an $R(U)$-covering of $A_t(U)$, such that if $g$ is a
loop in $A_t(U)$ with $g(0) = g(1) = Id$ and which corresponds to an element
$\gamma$ of $\Pi(U) = R(U)$, then it can be lifted to a curve $\hat{g}$ in
$\hat{A}_t'(U)$ such that $\hat{g}(0)$ is the unity element in $\hat{A}_t'(U)$ but
$\hat{g}(1)$ is the image of $\gamma$ in $\hat{A}_t'(U)$ via the injection $R(U) \to
\hat{A}_t'(U)$. The association $U \mapsto \hat{A}_t'(U)$ is a presheaf, which can
be made in to a sheaf $\hat{\mathcal A}_t$ whose stalk at each point $x \in O$ is
$\lim_{U(x) \to x} \hat{A}_t'(U(x))$. Clearly, $\hat{\mathcal A}_t$ is an extension
of ${\mathcal A}_t$ by ${\mathcal R}$: $0 \to {\mathcal R} \to \hat{\mathcal A}_t
\to {\mathcal A}_t \to 0 $.

We will define a subgroup $A_c(U)$ of $A_t(U)$, consisting of the
so-called {\it compatible} (with respect to the local ${\Bbb
S}^1$-actions) homeomorphisms of $(\pi^{-1}(U), {\mathcal L})$. For
this we will fix an appropriate thickened stratification $O = \sqcup
\bar{S_i}$ of the base space $O$ (recall that such a thickened
stratification is unique topologically), and assume that $U$ is
compatible with this thickened stratification, i.e. the intersection
of $\bar{S_i}$ with $U$ gives rise to a thickened stratification of
$U$. Let $y$ be a point in $U$. Then there is a point $x$ in $U$ such
that $y$ and $x$ belong to the same thickened stratum, and that the
stratum that contains $x$ has the same index as the thickened stratum
that contains $x$. Let $V(x)$ be a star-shaped neighborhood of $x$
which contains $y$ and which lies in the thickened stratum that
contains $x$. Then there is a natural system-preserving action of
${\Bbb T}^r$ on $\pi^{-1}(V(x))$, where $r = {\rm rank}\ R(V(x))$,
which is unique up to automorphisms of ${\Bbb T}^r$. In particular,
there is a natural action of the torus ${\Bbb T}^r$ on $\pi^{-1}(y)$
which is unique up to automorphisms of ${\Bbb T}^r$ and which does
not depend on the choice of $x$. We will denote this action on
$\pi^{-1}(y)$ (up to automorphisms) by $T(y)$. We will say that an
element $\phi$ of $A_t(U)$ is a {\it compatible homeomorphism} if it
commutes with the torus action $T(y)$ on $\pi^{-1}(y)$ for each $y \in U$.
Evidently, the set of compatible homeomorphisms is a subgroup of
$A_t(U)$, which we will denote by $A_c(U)$. The association
$U \mapsto A_c(U)$ gives rise to a non-Abelian subsheaf of the sheaf
${\mathcal A}_t$ over $O$, which we will denote by ${\mathcal A}_c$.
Of course, the definition of ${\mathcal A}_c$ depends on the choice
of a thickened stratification of $O$, but since such a stratification
is unique up to homeomorphisms, the sheaf ${\mathcal A}_c$ is also
unique up to isomorphisms.

Similarly, if $\pi : (\pi^{-1}(U(x)), \omega, {\mathcal L}) \to U(x)$
and $\pi : (\pi^{-1}(V(y)), \omega, {\mathcal L}) \to V(y)$ are two
topologically equivalent singularities, then we can talk about
compatible homeomorphisms from $(\pi^{-1}(U(x)), {\mathcal L})$
to $(\pi^{-1}(V(y)), {\mathcal L})$ (with respect to a given thickened
stratification of $U(x)$ and $V(y)$). Our next hypothesis is that such
compatible homeomorphisms always exist.

\begin{quote}
{\bf (H4)} {\it If two singularities
    $(\pi^{-1}(U(x)), \omega, {\mathcal L})$
and $(\pi^{-1}(V(y)), \omega, {\mathcal L})$
are topologically equivalent, then for any fibration-preserving
homeomorphism $\psi$ from $(\pi^{-1}(U(x)),{\mathcal L})$
to $(\pi^{-1}(V(y)), {\mathcal L})$ which maps
a given thickened stratification of $U(x)$ to a given thickened
 stratification of $V(y)$, then there is a compatible (with
 respect to the local ${\Bbb S}^1$-actions and given thickened stratifications)
homeomorphism $\phi$ from $(\pi^{-1}(U(x)),{\mathcal L})$
to $(\pi^{-1}(V(y)),{\mathcal L})$ such that $\phi^{-1} \circ \psi$
belongs to $A_t(U(x))$. }
\end{quote}

Hypothesis (H4) can be verified directly for nondegenerate singularities.
Similarly to the case of $A_t(U)$, the group $A_c(U)$ also has a natural
central extension by $R(U)$, which we will denote by $\hat{A}_c'(U)$. In
terms of sheaves, ${\mathcal A}_c$ has an extension by ${\mathcal R}$, which
we will denote by $\hat{\mathcal A}_c$, and we have the following commutative
diagram :
$$
\begin{array}{ccccccccc}
0 & \to & {\mathcal R} & \to & \hat{\mathcal A}_c & \to & {\mathcal A}_c & \to & 0 \\
  &     &   \parallel  &     & \downarrow         &     & \downarrow     &     &   \\
0 & \to & {\mathcal R} & \to & \hat{\mathcal A}_t & \to & {\mathcal A}_t & \to & 0
\end{array}
$$

Our next hypothesis is:

\begin{quote}
{\bf (H5)}{\it If $\pi : (\pi^{-1}(U(x)), \omega, {\mathcal L}) \to U(x)$
is a singularity then the group $\hat{A}_c'(U(x))$ is contractible.
Equivalently, $A_c(U(x))$ is homotopic to ${\Bbb T}^k$ where
$k = {\rm rank }\ R(U(x))$.
}
\end{quote}

Again, for nondegenerate singularities, Hypothesis (H5) can be verified
directly, and is similar to the fact that the space of homeomorphisms
from a cube to itself which are identity on the boundary is contractible
(with respect to the open-compact topology).

\subsection{Differential forms on the base space}
\label{subsection:DeRham}

Recall that a differential $k$-form $\beta$ on $(M^{2n}, {\mathcal L})$ is called a
{\it basic form} (with respect to the fibration ${\mathcal L}$) if for any vector
$X$ on $M^{2n}$ tangent to $\mathcal L$ we have $i_X \beta = i_X d\beta = 0$. We
will denote the space of basic $k$-forms on $M^{2n}$ (with respect to a given
integrable system $\pi : (M^{2n}, \omega, {\mathcal L}) \to O$) by $\Omega^k(O)$ and
consider it as the space of differential $k$-forms on the base space $O$. (If $O$ is
regular then basic forms on $M$ are pull-backs of differential forms on $O$). In
particular, the space of smooth functions on $O$, denoted by $C^{\infty}(O)$, is the
space of functions $f : O \to {\Bbb R}$ such that $f \circ \pi$ is a smooth function
on $M^{2n}$. For each $k$, the linear space $\Omega^k(O)$ is a
$C^{\infty}(O)$-module. The differential of a basic form is again a basic form. Thus
we have the following DeRham complex of $O$:
$$
0 \to C^{\infty}(O) \stackrel{d}{\to} \Omega^1(O) \stackrel{d}{\to} \Omega^2(O)
\stackrel{d}{\to} \cdots \stackrel{d}{\to} \Omega^n(O) \stackrel{d}{\to} 0
$$
The cohomology of this complex, known as {\it basic cohomology}, will also be called
the {\it DeRham cohomology of $O$}, and denoted by $H^k_{DR}(O, {\Bbb R})$. Of
particular interest is the second cohomology group $H^2_{DR}(O, {\Bbb R})$: if
$\beta$ is a closed 2-form on $O$ then $\omega + \pi^{\ast} \beta$ will be a new
symplectic form on $M$ for which $\pi : (M, \omega + \pi^{\ast}\beta, {\mathcal L})
\to O$ remains an integrable Hamiltonian system (see Subsection
\ref{subsection:magnetic}).

If all singularities over $O$ are nondegenerate, then it is easy to see that
$C^{\infty}$ is a fine sheaf over $O$, the Poincar\'e lemma holds for $O$ (i.e.
closed differential forms on $O$ are locally exact), and the DeRham cohomology of
$O$ is naturally isomorphic to its singular cohomology. Indeed, near each
nondegenerate singular point of the system in the symplectic manifold there is a
radiant-like vector field which preserves the fibration, and one can use these
radiant vector fields to prove the Poincar\'e lemma for $O$ by the same Thom-Moser
path method as used in the book of Abraham and Marsden \cite{AM}.

When dealing with the symplectic classification of general integrable Hamiltonian
systems, we will make the following hypothesis:

\begin{quote}
{\bf (H6)} {\it The sheaf $C^{\infty}$ of local smooth functions on $O$ is a fine
sheaf.}
\end{quote}

Notice that, when $O$ contains degenerate singularities, closed 1-forms on $O$ are
still locally exact: if $\beta$ is a closed 1-form in a star-shaped neighborhood
$U(x)$ of $O$, then since $\pi^{\ast}\beta$ vanishes on $\pi^{-1}(x)$ and
homotopically $\pi^{-1}(U(x))$ is the same as $\pi^{-1}(x)$, the cohomological class
of $\pi^{\ast}\beta$ in $H^1(\pi^{-1}(U(x)), {\Bbb R})$ is zero. Hence
$\pi^{\ast}\beta = d h$ for some function $h$ on $\pi^{-1}(U(x))$. It is obvious
that $h$ is constant on the fibers of the fibration, i.e. $h$ is a basic function.
As a consequence, we have the following short exact sequence of Abelian sheaves over
$O$:
$$
0 \longrightarrow {\Bbb R} \longrightarrow {C^{\infty}} \longrightarrow
{\mathcal Z}^1  \longrightarrow 0,
$$
where ${\mathcal Z}^1$ denotes the sheaf of local closed differential 1-forms.
Since $C^{\infty}$ is a fine sheaf by Hypothesis (H6), the cohomology long exact
sequence associated to the above short exact sequence gives rise to a natural
isomorphism from $H^k(O, {\mathcal Z}^1)$ to $H^{k+1}(O, {\Bbb R})$ ($k \geq 1$).

\subsection{Local symplectomorphisms}

Let $X$ be a  symplectic vector field on a singularity $(\pi^{-1}(U(x)), \omega,
{\mathcal L})$, which is tangent to  ${\mathcal L}$. Then the time-1 map of $X$,
denoted by $g_X^1$, is an element of the group $A_t(U(x))$ of admissible
homeomorphisms, which preserves the symplectic form. The inverse is also true, at
least for nondegenerate singularities:

\begin{lem}
Let $(\pi^{-1}(U(x)), \omega, {\mathcal L})$ be a nondegenerate singularity, and
$\phi$ a symplectomorphism from $(\pi^{-1}(U(x)), \omega, {\mathcal L})$ onto itself
which preserves every stratum of every fiber of the system. Then $\phi$ can be
written as the time one map $g_X^1$ of a symplectic vector field $X$ on
$(\pi^{-1}(U(x)), \omega)$ which is tangent to $\mathcal L$.
\end{lem}

{\it Sketch of the proof}. If the singularity contains only elliptic
components, then one can write an explicit formula for $X$ in a
canonical system of coordinates. If the singularity does not contain
any elliptic component, then the proof is similar to that of
Lemma \ref{lem:hyperbolic}. The general case is a combination of these
two cases. $\blacksquare$

We will make the above property into a hypothesis, our last one, which
will be used in the symplectic classification of integrable Hamiltonian
systems:

\begin{quote}
{\bf (H7)} {\it Let $(\pi^{-1}(U(x)), \omega, {\mathcal L})$ be a singularity, and
$\phi$ a symplectomorphism from $(\pi^{-1}(U(x)), \omega, {\mathcal L})$ onto itself
which preserves every stratum of every fiber of the system. Then $\phi$ can be
written as the time one map $g_X^1$ of a symplectic vector field $X$ on
$(\pi^{-1}(U(x)), \omega)$ which is tangent to $\mathcal L$ }.
\end{quote}

For each open subset $U$ of $O$, we will denote by $A_s(U)$ the set of smooth
symplectomorphisms from $(\pi^{-1}(U), \omega)$ onto itself which belong to the
group of compatible homeomorphisms $A_t(U)$. It is obvious that $A_s(U)$ is a
subgroup of the group of compatible homeomorphisms $A_c(U)$ (for any thickened
stratification of $O$), and the association $U \mapsto A_s(U)$ is a sheaf over $O$,
which we will call {\it the sheaf of local symplectomorphisms of the system}, and
denote by ${\mathcal A}_s$.

There is a natural projection from the sheaf ${\mathcal Z}^1$ of local
closed 1-forms on $O$ to ${\mathcal A}_s$: if $\beta$ is a basic
closed 1-form on $(\pi^{-1}(U(x)), \omega, {\mathcal L})$, then the
vector field $X$ defined by $i_X \omega = \beta$ is symplectic and its
time-1 map defines an element in $A_s(U(x))$. On the other hand, there
is a natural injection from the sheaf ${\mathcal R}$ of local symplectic
${\Bbb S}^1$-actions to ${\mathcal Z}^1$: if $g$ is a symplectic
${\Bbb S}^1$-action on $(\pi^{-1}(U(x)), \omega, {\mathcal L})$
and $X$ the symplectic vector field that generates $g$, then $i_X\omega$
is an element of ${\mathcal Z}^1(U(x))$. Hypothesis (H7) implies that
the projection from ${\mathcal Z}^1$ to ${\mathcal A}_s$ is surjective,
and its kernel is the image of $\mathcal R$ in ${\mathcal Z}^1$. In other
words, we have the following exact sequence of Abelian sheaves
$$
0 \longrightarrow {\mathcal R} \longrightarrow {\mathcal Z}^1
\longrightarrow {\mathcal Z}^1/{\mathcal R} \longrightarrow 0,
$$
and ${\mathcal A}_s$ is naturally isomorphic to
${\mathcal Z}^1/{\mathcal R}$.
\\


\section{Characteristic classes and classification}
\label{section:classification}

\subsection{Homological monodromy and rough equivalence}
\label{subsection:rough}

To define the characteristic classes of an integrable system, we will have to
compare it to another ``reference" system which is ``roughly equivalent'' to it.
Here two integrable Hamiltonian systems will be called roughly equivalent if they
have the same base space, the same singularities, and the same ``homological
monodromy'', in the following sense:

\begin{define}
Two integrable Hamiltonian systems $(M_a, \omega_a, {\mathcal
L}_a) \stackrel{\pi_a}{\longrightarrow} O_a$ and $(M_b, \omega_b,
{\mathcal L}_b) \stackrel{\pi_b}{\longrightarrow} O_b$ are called
{\it roughly topologically equivalent} if there is a homeomorphism
$\phi: O_a \to O_b$, a covering of $O_a$ by open subsets $U_i$, a
homeomorphism $\Phi_i: \pi_a^{-1}(U_i) \to \pi_b^{-1}(\phi (U_i))$
for each $i$, such that $\pi_b \circ \Phi_i = \phi \circ
\pi_a|_{\pi_a^{-1}(U_i)}$, and $\Phi_i^{-1} \Phi_j$ induces the
identity map on the fundamental groups of the strata of
$\pi_a^{-1}(x)$ and the identity map on $H_1(\pi_a^{-1}(x), {\Bbb
Z})$ for each point $x \in U_i \cap U_j$. The two systems are
called {\it roughly symplectically equivalent} if, in addition,
$\Phi_i$ are smooth symplectomorphisms.
\end{define}

It follows from the above definition and the hypotheses made in the previous section
that, if two systems are roughly topologically equivalent, then they have the same
affine monodromy sheaf $\mathcal R$. For regular Lagrangian torus fibrations, the
condition of rough topological equivalence is equivalent to the condition of having
the same base space up to homeomorphisms and the same affine monodromy.

If two systems are roughly symplectically equivalent, then their base spaces are
diffeomorphic and have the same integral affine structure, in the sense that there
is a homeomorphism between them which preserves the algebra of differential forms
and the sheaf of local integral affine functions. For two roughly topologically
equivalent systems to be roughly symplectically equivalent, a necessary and
sufficient condition is that their base spaces have the same affine structure and
their singularities are not only topologically equivalent but also symplectically
equivalent. Symplectic invariants of simplest singularities and systems (with one or
one and a half degrees of freedom) have been studied by some authors (see e.g.
\cite{BF,DMT,Vu3}).

If two systems are roughly topologically (resp. symplectically)
equivalent, then we will also say that they have the same {\it rough
topological (resp. symplectic) type}. Given an integrable Hamiltonian
system $ \pi : (M, \omega, {\mathcal L}) \to O $, we will denote its
rough topological and symplectic types by $\hat{O}_{top}$ and
$\hat{O}_{symp}$ respectively, and view them as {\it framed base spaces},
with the framing given by singularities and ``monodromies''. For each rough
topological or symplectic type, we will try and choose a ``reference
system" with this rough type. For example, if the base space $O$ is
regular, then the framing is given by the affine monodromy,
and the obvious reference system is the one which admits a
global (Lagrangian) section. In the general case, where there is no
obvious choice, we'll pick an arbitrary system and call it the
reference system for a given rough type.

\subsection{Definition of characteristic classes}
\label{subsection:definitions}

Assume that an integrable Hamiltonian system $\pi : (M, \omega,
{\mathcal L}) \to O$ is roughly topologically equivalent to a
reference system $\pi_0 : (M_0, \omega_0, {\mathcal L}_0) \to O$ over
the same base space $O$. By definition, there is a covering of $O$ by
open subsets $U_i$, a homeomorphism $\Phi_i: \pi_0^{-1}(U_i) \to
\pi^{-1}(\phi (U_i))$ for each $i$, such that $\pi \circ \Phi_i =
\pi_0|_{\pi_0^{-1}(U_i)}$, and $\mu_{ij} = \Phi_i^{-1} \circ \Phi_j $
belongs to the group $A_t(U_i \cap U_j)$ (defined in the previous section,
for the reference system). Its clear that the family
$(\mu_{ij} \in A_t (U_i \cap U_j))$ is a Cech 1-cocycle in the
non-Abelian sheaf ${\mathcal A}_t$, which will define an element $\mu$
in $H^1({\mathcal A}_t)$. The first cohomology class $\mu$ is a
topological invariant of the system $\pi : (M, \omega,
{\mathcal L}) \to O$, and it is trivial if and only if
$\pi : (M, \omega, {\mathcal L}) \to O$ is topologically equivalent
to the reference system $\pi_0 : (M_0, \omega_0, {\mathcal L}_0) \to O$
by a homeomorphism which projects to the identity map on $O$.

Similarly, if $\pi : (M, \omega, {\mathcal L}) \to O$ is roughly
symplectically equivalent to a given reference system $\pi_0 : (M_0,
\omega_0, {\mathcal L}_0) \to O$, then we will be able to define
an element $\mu_L$ in $H^1({\mathcal A}_s)$, where ${\mathcal A}_s$
is the sheaf of local symplectomorphisms of the reference system.
Recall that ${\mathcal A}_s \cong {\mathcal Z}^1/{\mathcal R}$
(and is an Abelian sheaf). Thus we can write
$\mu_L \in H^1({\mathcal Z}^1/{\mathcal R})$.


For a given reference system, the short exact sequence
$ 0  \to  {\mathcal R}  \to  \hat{\mathcal A}_{t}  \to  {\mathcal A}_{t}  \to  0 $
gives rise to the following long exact sequence:
$$
...  \to  H^1(\hat{\mathcal A}_{t})  \to
H^1({\mathcal A}_{t})  \stackrel{\delta}{\to}  H^2({\mathcal R})  \to
H^2(\hat{\mathcal A}_{t})  \to ...
$$

\begin{define}
The second cohomology class $\mu_C \in H^2(O, {\mathcal R})$ which is the image
of $\mu$ under the coboundary map $\delta: H^1(O, {\mathcal A}_t) \to
H^2(O, {\mathcal R})$ arising from the short exact sequence
$0 \to {\mathcal R} \to \hat{\mathcal A}_{t} \to {\mathcal A}_{t} \to  0 $
is called the {\it Chern class} of the system $(M, \omega, {\mathcal L})$
(with respect to a given reference system $(M_0, \omega_0, {\mathcal L}_0)$).
The first cohomology class $\mu_{L} \in H^1(O, {\mathcal Z}^1/{\mathcal R})$
is called the {\it Lagrangian class} of $(M, \omega, {\mathcal L})$ .
\end{define}

It is clear from the above definition that the Chern class is a topological
invariant: if two systems are topologically equivalent, then they will have
the same Chern class (with respect to any common reference system) after an
appropriate homeomorphism of their respective base spaces. Similarly,
the Lagrangian class is obviously a symplectic invariant. Notice that
in the definition of the characteristic classes, we need a reference system,
therefore the Chern and Lagrangian classes are not ``absolute'' classes but
rather ``relative" ones. In other words, they live in an affine space rather
than a vector space.

\subsection{Chern class via compatible homeomorphisms}

It follows from Hypothesis (H4) that, in the definition of the
Chern class, we may replace the sheaf ${\mathcal A}_t$ of local
admissible homeomorphisms by the sheaf ${\mathcal A}_c$ of local
compatible (with respect to a given thickened stratification of
the base space) homeomorphisms: we can choose homeomorphisms
$\Phi_i: \pi_0^{-1}(U_i) \to \pi^{-1}(\phi (U_i))$ such that
$\mu_{ij} = \Phi_i^{-1} \circ \Phi_j \in A_c(U_i \cap U_j)$. Thus
the cocycle  $\mu_{ij}$ defines an element $\bar{\mu} \in H^1(O,
{\mathcal A}_c)$ which maps to $\mu \in H^1(O, {\mathcal A}_t)$
under the natural homomorphism $H^{\ast}(O, {\mathcal A}_c) \to
H^{\ast}(O, {\mathcal A}_t)$. Under the long exact sequence
$$
...  \to  H^1(\hat{\mathcal A}_{c})  \to
H^1({\mathcal A}_{c})  \stackrel{\bar{\delta}}{\to}  H^2({\mathcal R})  \to
H^2(\hat{\mathcal A}_{c})  \to  ...
$$
we have $\bar{\delta}\bar{\mu} = \mu_C$.

If $\bar{\mu} = 0$ then of course the system $(M,\omega, {\mathcal L})$ is
topologically equivalent to the reference system $(M_0,\omega_0, {\mathcal L}_0)$,
and vice versa. If we can show that $H^1(O, \hat{\mathcal A}_c) = 0$, then we will
have $\bar{\mu} = 0$ if and only if the Chern class $\mu_C$ vanishes.

\begin{lem}
\label{lem:acyclic}
The sheaf $\hat{\mathcal A}_c$ is acyclic. In particular, we have
$H^1(O, \hat{\mathcal A}_c) = 0$.
\end{lem}

{\it Proof}. The idea is that, due to Hypothesis (H5), the sheaf
$\hat{\mathcal A}_c$ is similar to a fine sheaf, though it is not
Abelian. It suffices to prove the above lemma for an arbitrary
thickened stratum in $O$ (more precisely, an arbitrary thickened
open subset of a stratum). The acyclicity of $\hat{\mathcal A}_c$
over $O$ will then follow from the Meyer-Vietoris exact sequences.
Let $\bar{S}$ be a thickened stratum of $O$, and consider the
sheaf $\hat{\mathcal A}_c$ over $\bar{S}$. $\bar{S}$ may be
considered as a locally trivial fibration over a stratum $S$, with the fiber
being local stratified submanifolds in $\bar{S}$ which are transversal
to $S$. We will consider only open subsets of $\bar{S}$ which
are saturated by the fibers of this fibration. In other words, we will
consider $\hat{\mathcal A}_c$ as a sheaf over the $S$. Then it will
become the sheaf of local continuous sections of a locally trivial fibration
over $S$ whose fibers are isomorphic to a contractible topological group,
which implies the acyclicity. $\blacksquare$

\subsection{Additivity of characteristic classes}

The Chern class is additive in the following sense:

\begin{lem}
\label{lem:additive}
If three systems $\pi_i : (M_i,\omega_i,{\mathcal L}_i) \to O$, $i= 1,2,3$,
over the same base space $O$, are roughly topologically equivalent, and the
Chern class of $(M_i,\omega_i,{\mathcal L}_i)$ with respect to
$(M_j,\omega_j,{\mathcal L}_j)$ is
$\mu^{ij}_C \in H^2(O, {\mathcal R})$ (the three systems have the same
affine monodromy sheaf $\mathcal R$), then we have
$\mu^{12}_C + \mu^{23}_C + \mu^{31}_C = 0$.
\end{lem}

{\it Proof}. The proof is straightforward and uses the fact that
${\mathcal R}$ lies in the center of $\hat{\mathcal A}_t$: Denote
by $(U_i)$ an appropriate open covering of $O$, $U_{ij} = U_i \cap
U_j$, ${\phi}_i^{ab} : \pi_a^{-1}(U_i) \to \pi_b^{-1}(U_i)$ \
($a,b= 1,2,3$) the homeomorphisms that define the rough
topological equivalences. We can assume that ${\phi}_i^{ba} =
({\phi}_i^{ab})^{-1} $ and ${\phi}_i^{31} \circ {\phi}_i^{23}
\circ {\phi}_i^{12} = Id$. Put ${\phi}_{ij}^{ab} =
({\phi}_i^{ab})^{-1} \circ {\phi}_j^{ab} \in A^a_t(U_{ij})$, and
denote by ${\psi}_{ij}^{ab}$ a lifting of ${\phi}_{ij}^{ab}$ from
the group $A^a_t(U_{ij})$ to the group $\hat{A}^a_t(U_{ij})$ such
that ${\psi}_{ji}^{ab} = ({\psi}_{ij}^{ab})^{-1}$. Then
$\mu_{ijk}^{ab} = {\psi}_{ij}^{ab} {\psi}_{jk}^{ab}
{\psi}_{ki}^{ab} \in R(U_{ijk})$ is a 2-cocycle that defines the
Chern class of the system $(M_b,\omega_b,{\mathcal L}_b)$ with
respect to the system $(M_a,\omega_a,{\mathcal L}_a)$. If we can
choose the elements ${\psi}_{ij}^{ab}$ in such a way that
$\mu_{ijk}^{12} + \mu_{ijk}^{23} + \mu_{ijk}^{31} = 0$ then we are
done. Since ${\phi}_{ij}^{23} = {\phi}_{i}^{12} \circ
{\phi}_{ij}^{13} ({\phi}_{ij}^{12})^{-1} \circ
({\phi}_{i}^{12})^{-1}$, we can choose ${\psi}_{ij}^{ab}$ so that
${\psi}_{ij}^{23} = Ad_{{\phi}_{i}^{12}} ({\psi}_{ij}^{13}
({\psi}_{ij}^{12})^{-1})$ (the operator $Ad_{{\phi}_{i}^{12}} :
\hat{A}^1_t(U_{ij}) \to \hat{A}^2_t(U_{ij}))$ is well defined).
Then we have: $ \mu_{ijk}^{23} = {\psi}_{ij}^{23} {\psi}_{jk}^{23}
{\psi}_{ki}^{23} = Ad_{{\phi}_{i}^{21}} ({\psi}_{ij}^{23}
{\psi}_{jk}^{23} {\psi}_{ki}^{23}) $ (because
$Ad_{{\phi}_{i}^{21}}$ is identity when restricted to
$R(U_{ijk})$) $  = ({\psi}_{ij}^{13} ({\psi}_{ij}^{12})^{-1})
\circ Ad_{{\psi}_{ij}^{12}} ({\psi}_{jk}^{13}
({\psi}_{jk}^{12})^{-1}) \circ Ad_{{\psi}_{ik}^{12}}
({\psi}_{ki}^{13} ({\psi}_{ki}^{12})^{-1}) $ $ = {\psi}_{ij}^{13}
{\psi}_{jk}^{13} ({\psi}_{jk}^{12})^{-1} ({\psi}_{ij}^{12})^{-1}
{\psi}_{ik}^{12} {\psi}_{ki}^{13} $ $ =  ({\psi}_{ki}^{13}
{\psi}_{ij}^{13} {\psi}_{jk}^{13}) (({\psi}_{jk}^{12})^{-1}
({\psi}_{ij}^{12})^{-1} {\psi}_{ik}^{12}) $ $ = \mu_{ijk}^{13} -
\mu_{ijk}^{12}. $ $\blacksquare$

Similarly, the Lagrangian class is also additive:

\begin{lem}
\label{lem:additive2}
If three systems $\pi_i : (M_i,\omega_i,{\mathcal L}_i) \to O$, $i= 1,2,3$,
over the same base space $O$, are roughly symplectically equivalent, and the
Lagrangian class of $(M_i,\omega_i,{\mathcal L}_i)$ with respect to
$(M_j,\omega_j,{\mathcal L}_j)$ is $\mu^{ij}_L \in H^1(O, {\mathcal Z}^1/{\mathcal R})$
(the three systems have the same sheaf ${\mathcal Z}^1/{\mathcal R}$), then we have
$\mu^{12}_L + \mu^{23}_L + \mu^{31}_L = 0$.
\end{lem}

The proof is obvious. $\blacksquare$

\subsection{Classification theorems}
\label{subsection:classification}

We can now present the main theorems. Their formulations are similar to the case
of regular Lagrangian torus fibrations.

\begin{thm}
If two roughly topologically equivalent integrable Hamiltonian systems have the
same Chern class (with respect to a common reference system), then they
are topologically equivalent.
\end{thm}

The proof of the above theorem follows directly from Lemma \ref{lem:acyclic}
and Lemma \ref{lem:additive}. $\blacksquare$

\begin{thm}
If two roughly symplectically equivalent integrable Hamiltonian systems have the
same Lagrangian class (with respect to a common reference system), then they
are symplectically equivalent.
\end{thm}

The proof of the above theorem follows directly from the definition. $\blacksquare$

For a given reference system, the following commutative diagram of short
exact sequences,
$$ \begin{array}{ccccccccc}
0 & \to & {\mathcal R} & \to & \hat{\mathcal A}_{t} & \to & {\mathcal A}_{t} & \to & 0 \\
  &     &   ||     &     & \uparrow           &     & \uparrow       &  \\
0 & \to & {\mathcal R} & \to & {\mathcal Z}^1         & \to &{\mathcal A}_{s} & \to & 0
\end{array}
$$
give rise to the following commutative diagram of associated long exact sequences
of cohomologies over the base space $O$:
$$
\begin{array}{ccccccccccc}
... & \to & H^1(\hat{\mathcal A}_{t}) & \to &
H^1({\mathcal A}_{t}) & \stackrel{\delta}{\to} & H^2({\mathcal R}) & \to &
H^2(\hat{\mathcal A}_{t}) & \to & ... \\
    &     & \uparrow & & \uparrow & & || & & \uparrow & & \\
... & \stackrel{\hat{d}}{\to} & H^1({\mathcal Z}^1) & \to &
H^1({\mathcal A}_{s}) & \stackrel{\Delta}{\to} & H^2({\mathcal R}) &
\stackrel{\hat{d}}{\to} & H^2({\mathcal Z}^1) & \to & ... \\
& &   ||               & & & & & &      ||            & & \\
& & H^2(O, {\Bbb R}) & & & & & & H^3(O, {\Bbb R}) & &
\end{array}
$$
In the above diagram, $H^k(O,{\mathcal Z}^1)$ are identified with
$H^{k+1}(O,{\Bbb R})$ via the isomorphisms arising from the short exact
sequence $0 \to {\Bbb R} \to C^{\infty} \stackrel{d}{\to} {\mathcal Z}^1 \to 0$.
Remark that, in the above diagram, the operators
$\delta$ depend on the topological type of the reference system
$\pi_0 : (M_0, \omega_0, {\mathcal L}_0) \to O$, while
the operators $\hat{d}$ and $\Delta$ depend on its rough symplectic type.

It follows from the above commutative diagram that if the system
$(M, \omega, {\mathcal L})$ is roughly symplectically equivalent to a reference
system $(M_0, \omega_0, {\mathcal L}_0)$, and both the Chern and the Lagrangian
classes of $(M, \omega, {\mathcal L})$ are taken with respect to
$(M_0, \omega_0, {\mathcal L}_0)$, then under the maps $\Delta$ and
$\hat{d}$ we have $\mu_{L} \stackrel{\Delta}{\rightarrow} \mu_{C}
\stackrel{\hat{d}}{\rightarrow} 0$. It is clear from the construction of
characteristic classes that any element in $H^1(O, {\mathcal Z}^1/{\mathcal R})$ is
the Lagrangian class of some integrable Hamiltonian system which
is roughly symplectically equivalent to $(M_0, \omega_0, {\mathcal L}_0)$.
Thus we have the following proposition which is similar to a result of
Dazord and Delzant \cite{DD} for the regular case:

\begin{prop}
An element $\mu_C \in H^2(O, {\mathcal R})$
is the Chern class of some integrable Hamiltonian system
$\pi :(M, \omega, {\mathcal L}) \to O$ roughly
symplectically equivalent to a given reference system
$\pi_0 :(M_0, \omega_0, {\mathcal L}_0) \to O$ over a given base
space $O$ if and only if $\hat{d}(\mu_C) = 0$.
Under this condition, the space of integrable
Hamiltonian systems which are roughly symplectically equivalent to
$\pi_0 : (M_0, \omega_0, {\mathcal L}_0) \to O$
and which have $\mu_C$ as the Chern class, considered up to symplectic
equivalence, is naturally isomorphic to $H^2(O, {\Bbb R}) /
\hat{d}(H^1(O,{\mathcal R}))$.
\end{prop}

Notice that the above proposition does not solve the following problem:
which elements in $H^2(O, {\mathcal R})$ can be realized as the Chern class
of a system which is roughly topologically equivalent to
$\pi_0 :(M_0, \omega_0, {\mathcal L}_0) \to O$ ? The condition
$\hat{d}(\mu_C) = 0$ needs not hold for systems which are roughly
topologically equivalent but not roughly symplectically equivalent
to the reference system (except for the case $H^3(O, {\Bbb R}) = 0$,
when this condition is empty). I don't know which conditions must
$\mu_C$ satisfy in general.

\subsection{The magnetic term}
\label{subsection:magnetic}

In analogy with classical electromagnetism (see e.g. \cite{GS2}), by a
{\it magnetic term} we will mean (the pull-back $\pi^{\ast}\beta$ of) a
closed 2-form $\beta$ on the base space $O$ of an integrable system
$\pi : (M, \omega, {\mathcal L}) \to O$. We have:

\begin{lem}
For any closed 2-form $\beta$ on $O$, the form $\omega + \pi^{\ast}\beta$
is a symplectic form on $M$, and
$\pi : (M, \omega + \pi^{\ast}\beta, {\mathcal L}) \to O$ is an integrable system.
\end{lem}

{\it Proof}. We will show that $\omega + \pi^{\ast}\beta$ is nondegenerate
everywhere on $M$, the rest is obvious. Let $p$ be an arbitrary point of $M$.
We will show that the kernel $K_p$ of $\pi^{\ast}\beta$ at $p$ is a coisotropic
subspace of the tangent space at $p$ with respect to the symplectic form
$\omega$. Then for any vector $X \in K_p$ there is a vector $Y \in T_pM$
such that $\omega(X,Y) + \pi^{\ast}\beta (X,Y) = \omega(X,Y) \neq 0$, and
for any vector $X \in T_pM$ not belonging to $K_p$ there is a vector
$Y \in K_p$ such that $\omega(X,Y) + \pi^{\ast}\beta (X,Y) = \omega(X,Y) \neq 0$.
According to Hypotheses (H1) and (H2), there is a sequence  of
regular points $p_i \in M$, $i \in {\Bbb N}$, of the fibration ${\mathcal L}$ which
tend to $p$. Denote the tangent space of $\mathcal L$ at $p_i$ by $K_i$.
By taking a subsequence of $(y_p)$, we can assume that $K_i$ tends to a subspace
$K_0 \subset T_pM$ at $p$. Since each $K_i$ is Lagrangian, $K_0$ is also
Lagrangian. On the other hand, $K_i$ lies in the kernel of $\pi^{\ast}\beta$
at $p_i$, and it follows from the semi-continuity of the kernel that
$K_0 \subset K_p$. Thus $K_p$ is coisotropic. $\blacksquare$

It is easy to see that if the magnetic term $\beta$ is locally exact on  $O$,
then the two systems $\pi : (M, \omega + \pi^{\ast}\beta, {\mathcal L}) \to O$
and $\pi : (M, \omega, {\mathcal L}) \to O$ are roughly symplectically
equivalent, and if $\beta$ is exact then the two systems are symplectically
equivalent.

Assume in this paragraph that the Poincar\'e lemma holds for 2-forms on $O$,
i.e. any closed 2-form on $O$ is exact, and denote by ${\mathcal Z}^2$ the sheaf of local
closed 2-forms. Then we have the following short exact sequence:
$ 0  \to  {\mathcal Z}^1  \to  \Omega^1  \stackrel{d}{\to}  {\mathcal Z}^2  \to  0 $,
where the sheaf $\Omega^1$ of local 1-forms on $O$ is a fine sheaf
(because $C^{\infty}$ is a fine sheaf by Hypothesis (H6)). The associated
long exact sequence gives rise to the following isomorphisms:
$H^2_{DR}(O, {\Bbb R}) = H^0(O, {\mathcal Z}^2) / d H^0 (O, \Omega^1)
= H^1 (O, {\mathcal Z}^1) = H^2 (O, {\Bbb R})$. In other words, the second
DeRham cohomology group of $O$ is isomorphic to the second singular
cohomology group of $O$. This isomorphism, together with the natural operator
from $H^2 (O, {\Bbb R}) = H^1(O, {\mathcal Z}^1)$ to
$H^1(O, {\mathcal Z}^1/{\mathcal R})$, will send each closed 2-form
$\beta$ on $O$ to an element $\mu_L(\beta) \in H^1(O, {\mathcal Z}^1/{\mathcal R})$
which is nothing else but the Lagrangian class of the system
$\pi : (M, \omega + \pi^{\ast}\beta, {\mathcal L}) \to O$
with respect to the system $\pi : (M, \omega, {\mathcal L}) \to O$.
In particular, if the Poincar\'e lemma holds for 2-forms on $O$, then
any two systems over $O$ which are roughly symplectically equivalent
and have the same Chern class will differ from each other by only a magnetic term.

\subsection{Realization problem and integrable surgery}

Given a stratified manifold $O$ equipped with an integral affine structure
(in the sense of Subsection \ref{subsection:affine}), one can ask whether
it can be realized as the base space of some integrable Hamiltonian system.
If it is the case, we will say that $O$ is {\it realizable}. Of course,
if $O$ is to be realizable, it has to be {\it locally realizable}: each
singular point $y$ in $O$ corresponds to some singularity of some integrable
system, that is a singular Lagrangian torus fibration with the base space
$U(y)$ where $U(y)$ is a neighborhood of $y$ in $O$, in such a way that the
following compatibility is satisfied: If $U(y_1) \cap U(y_2) \neq \emptyset$
then there is a fibration-preserving symplectomorphism
$\Phi_{y_1y_2}$ between the restrictions of the two fibrations to
$U(y_1) \cap U(y2)$; if $U(y_1) \cap U(y_2) \cap U(y_3) \neq \emptyset$,
then the map $\Phi_{y_1y_2} \circ \Phi_{y_2y_3} \circ \Phi_{y_3y_1}$
(on a restricted fibration over $U(y_1) \cap U(y_2) \cap U(y_3)$) is
isomorphic to identity. A stratified integral affine manifold $O$
equipped with such singularities will be called a {\it formal rough
symplectic type}, and denoted by $\hat{O}_{symp}$ as before.

A natural problem arises: given a formal rough symplectic type $\hat{O}_{symp}$,
is there any integrable system roughly symplectically equivalent to it ?
A natural way to solve this problem is via {\it integrable surgery}.
In this paper, by an integrable surgery, we mean a surgery of an
integrable Hamiltonian system which projects to a surgery on its base
space.

For example, given two integrable systems
$(M_1,\omega_1, {\mathcal L}_1) \stackrel{\pi_1}{\to} O_1$ and
$(M_2,\omega_2, {\mathcal L}_2) \stackrel{\pi_2}{\to} O_2$
over two subsets $O_1$ and $O_2$ of a space $O$, such that they are
roughly symplectically equivalent when restricted to the common base space
$O_1 \cap O_2$. Does there exist an integrable system over
$O_1 \cup O_2$ which is roughly symplectically equivalent to the above
systems when restricted to $O_1$ and $O_2$ ? The answer to this question
may be given in terms of characteristic classes:

\begin{prop}
\label{prop:glue}
Denote the difference between the Lagrangian classes of the systems
$(M_1,\omega_1, {\mathcal L}_1) \stackrel{\pi_1}{\longrightarrow} O_1$ and
$(M_2,\omega_2, {\mathcal L}_2) \stackrel{\pi_2}{\longrightarrow} O_2$
restricted to $Q_1 \cap O_2$ by
$\mu_L \in H^1 (O_1 \cap O_2, {\mathcal Z}^1/{\mathcal R})$. Then there is an
integrable system with the base space $O_1 \cup O_2$ which is roughly symplectically
equivalent to the above two systems when restricted to $O_1$ and $O_2$ if and
only if $\mu_L$ lies in the sum of the images of
$H^1(O_1,{\mathcal Z}^1/{\mathcal R})$ and $H^1(O_2,{\mathcal Z}^1/{\mathcal R})$
in $H^1(O_1 \cap O_2,{\mathcal Z}^1/{\mathcal R})$ under the restriction maps.
In particular, if $(M_1,\omega_1, {\mathcal L}_1) \stackrel{\pi_1}{\longrightarrow} O_1$
and $(M_2,\omega_2, {\mathcal L}_2) \stackrel{\pi_2}{\longrightarrow} O_2$
are topologically equivalent when restricted to $O_1 \cap O_2$, then the obstruction
to the existence of the required system over $O_1 \cup O_2$ lies in the
quotient of the group
$H^2 (O_1 \cap O_2, {\Bbb R})/ \hat{d} H^1 (O_1 \cap O_2, {\mathcal R})$ by the sum
of the images of $H^2 (O_1, {\Bbb R})/ \hat{d} H^1 (O_1, {\mathcal R})$
and $H^2 (O_2, {\Bbb R})/ \hat{d} H^1 (O_2, {\mathcal R})$ under the restriction maps.
\end{prop}

{\it Proof}. It is a direct consequence of the results of
Subsection \ref{subsection:classification}
$\blacksquare$

For the case of systems with 2 degrees of freedom we have:

\begin{prop}
\label{prop:realization2}
Any formal rough symplectic type $\hat{O}_{symp}$ with $O$ two-dimensional
is realizable.
\end{prop}

{\it Proof}. If $O$ is 2-dimensional, then we can always choose $O_1 \cap O_2$
in Proposition \ref{prop:glue} to be a tubular neighborhood of something
1-dimensional, so all the obstructions vanish.
$\blacksquare$

Starting with 3 degrees of freedom, there are formal rough symplectic types which
are not realizable, as the following example shows.

\begin{example}
A fake base space \label{example:fake}
\end{example}

Let $(S^2,\omega)$ be a symplectic 2-sphere, $f: S^2 \to {\Bbb R}$ a Morse function
with 2 maximal points of the same value ($=1$), 2 minimal points of the same value
$(=-1)$, two saddle points of different values $(= \pm 1/2)$, such that $f$ is
invariant under an involution of $S^2$ which preserves the symplectic form. Denote
the base space of this integrable system with one degree of freedom by $G = G_+ \cup
G_-$, where $G_+$ (resp., $G_-$) corresponds to the part of the sphere with $f \geq
0$ (resp., $f \leq 0$). $G$ is a tree with 5 edges: 2 upper, 2 lower, and one
middle. Denote by $\sigma$ the involution of $G$ which preserves $f$ and the lower
edges but interchanges the two upper edges (so $\sigma$ cannot be lifted to an
involution on $S^2$). Denote by $K^2$ the Klein bottle with a standard integral
affine structure. We have $\pi_1(K^2) = <a,b \ | \ abab^{-1} = 1 >$. Denote by
$\bar{K}$ the double covering of $K^2$ corresponding to the subgroup of $\pi_1(K^2)$
which is generated by $a^2$ and $b$ (so $\bar{K}$ is again a Klein bottle), and
denote the involution on $\bar{K}$ corresponding to that double covering also by
$\sigma$. Put $O = \bar{K} \times_{\sigma} G = (\bar{K} \times G)/{\Bbb Z}_2$, with
the induced integral affine structure from the direct product. We have $O = O_+ \cup
O_-$, with $O_- = \bar{K} \times_{\sigma} G_- = K^2 \times G_-$, and $O_+ = \bar{K}
\times_{\sigma} G_+$ a twisted product. The spaces $O_-$ and $O_+$ are base spaces
of integrable systems induced from the direct product of the subsystems over $G_-$
and $G_+$ with a system over $\hat{K}$. These two systems are roughly equivalent
over $O_0 = O_+ \cap O_- \approx K^2$, but they are not equivalent, hence they
cannot be glued together to obtain a system over $O$. More precisely, the affine
monodromy sheaf over $O_0 \approx K^2$ in $O$ is ${\mathcal R} \cong {\mathcal
R}_{K^2} \oplus {\Bbb Z}$ where ${\mathcal R}_{K^2}$ is the affine monodromy of
$K^2$ as an affine manifold itself; $H^2 (O_0, {\mathcal R}) \cong H^2 (K^2,
{\mathcal R}_{K^2}) \oplus H^2 (K^2, {\Bbb Z}) = H^2(K^2, {\mathcal R}_{K^2}) \oplus
{\Bbb Z}_2$, and we have a natural map to the second component: $H^2 (O_0, {\mathcal
R}) \stackrel{\rho}{\longrightarrow} {\Bbb Z}_2$. Any system over $O_-$ will have a
Chern class which when restricted to $O_0$ will map to $0$ under the map $\rho$, but
any system over $O_-$ will have a Chern class which when restricted to $O_0$ will
map to the nontrivial element of ${\Bbb Z}_2$. hence those systems can never be
glued together to form a system over $O$. In other words, $O$ is not realizable.\\

\subsection{Some examples of integrable surgery}
\label{subsection:surgery}

\begin{example}Exotic symplectic ${\Bbb R}^{2n}s$
\end{example}

Start with the following two integrable systems: The first one is given by the
moment map ${\bf F} = (F_1,..., F_n) = (\pi x_1^2 + \pi y_1^2,..., \pi x_n^2 + \pi
y_n^2)$ on the open ball of radius 1 of ${\Bbb R}^{2n}$ with coordinates $x_i, y_i$
and with the standard metric and symplectic structure (i.e. a harmonic oscillator;
here $\pi = 3,14159...$). On the base space $O_1$ of this system, the functions
$F_i$ are also integral affine coordinates of the induced affine structures outside
the singularities. Let $O_2$ be an open $n$-disk attached to $O_1$ in such a way
that $O_1 \cup O_2$ is diffeomorphic to $O_1$ {\it rel.} singularities of $O_1$, and
$O_1 \cap O_2$ is contractible. Extend the functions $F_1,...F_n$ from $O_1$ to
$O_2$ in such a way that $dF_1 \wedge ... \wedge dF_n \neq 0$ everywhere on $O_2$
and there is a point $y \in O_2$ with $F_1 (y) = ... = F_n (y) = 0$. $O_2$ with the
integral affine structure given by the functions $F_i$ is the base space of a unique
integrable system $(M_2, \omega_2) \to O_2$. This is our second system. By
construction, our two systems can be glued in a unique natural way into an
integrable system living on a symplectic manifold diffeomorphic to ${\Bbb R}^{2n}$.
The preimage of $y$ in this manifold is a Lagrangian torus, and in fact it is an
{\it exact} Lagrangian torus (i.e. for any 1-form $\alpha$ such that $d\alpha$ is
equal to the symplectic form, the restriction of $\alpha$ on this torus is
cohomologous to 0). On the other hand, a famous result of Gromov \cite{Gromov} says
that in the standard symplectic space there can be no smooth closed exact Lagrangian
submanifold. Thus our symplectic space is exotic in the sense that it can not be
symplectically embedded to the standard symplectic space of the same dimension. This
example is inspired by a different example found by Bates and Peschke \cite{BP}.

\begin{example}Toric manifolds
\end{example}

Consider a Hamiltonian ${\Bbb T}^n$ action on a closed symplectic $2n$-dimensional
manifold $(M, \omega)$, which is free somewhere. $(M, \omega)$ together with this
torus action may be called a {\it Hamiltonian toric manifold}. The regular
(singular) orbits of this ${\Bbb T}^n$ action are Lagrangian (isotropic) tori, and
they are fibers of an integrable Hamiltonian system with only elliptic
singularities. The base space of this system is integral-affinely equivalent to a
polytope in the Euclidean space ${\Bbb R}^n$, whose each vertex has exactly $n$
edges and these edges can be moved to the principal axis of ${\Bbb R}^n$ by an
integral affine transformation. (This fact follows easily from the normal form of
elliptic singularities given by Eliasson \cite{Eliasson} and Dufour and Molino
\cite{DM}). A famous theorem of Delzant \cite{Delzant} says that each polytope
satisfying the above condition on vertices is the base space of a Hamiltonian toric
manifold which is unique up to symplectic equivalence. (These Hamiltonian toric
manifolds admit a K\"ahler structure which make them toric manifolds in the sense of
complex algebraic geometry). The uniqueness in Delzant's theorem is evident from our
point of view: Since ${\mathcal R}$ in this case is isomorphic to the constant sheaf
${\Bbb Z}^n$, and the base space is contractible, there is no room for
characteristic classes. The existence is also simple: one starts from a Lagrangian
section, and reconstructs the system (and the ambient manifold) in a unique way (see
\cite{BM}).

\begin{example}
Twisted products
\end{example}

We may call a {\it twisted product} of two integrable systems $(M_1,\omega_1,
{\mathcal L}_1) \stackrel{\pi_1}{\longrightarrow} O_1$ and $(M_2,\omega_2, {\mathcal
L}_2) \stackrel{\pi_2}{\longrightarrow} O_2$ an integrable system over $O_1 \times
O_2$, which is not topologically equivalent but is roughly symplectically equivalent
to the direct product of the two systems, and with the following property: The
Marsden-Weinstein reduction of this system to $\{y\} \times O_2$ (resp., $O_1 \times
\{y\}$) is symplectically equivalent to $(M_1,\omega_1, {\mathcal L}_1)$ (resp.,
$(M_2,\omega_2, {\mathcal L}_2)$ for every point $y \in O_1$ (resp. , $y \in O_2$).
For example, if $M_i$ if $M_i$ are symplectic tori, with the systems given by Morse
functions, then $H^2(O_1 \times O_2, {\mathcal R}) = {\Bbb Z}^4$ (here $\mathcal R$
is the corresponding affine monodromy sheaf, and the formula can be obtained easily
via Meyer-Vietoris sequences), and non-zero elements of this group corresponds to
twisted products.

\begin{example}
Blow-ups
\end{example}

Blowing-up, one of the main tools in algebraic geometry, is also
useful in symplectic geometry (see e.g. \cite{GS,McDuff,MS}), as
well as in the study of singularities of integrable Hamiltonian
systems (see e.g. \cite{DM}). To blow up a point in a $2n$-dimensional
symplectic manifold, one can cut away a symplectic ball containing this point,
and then collapse the sphere which is the boundary of this ball
to ${\Bbb C}P^{n-1}$ by collapsing each characteristic curve on this
sphere to a point. Since a symplectic ball admits a simple natural
integrable system, namely the harmonic oscillator, blowing-up can also be
done by integrable surgery: it amounts to cutting out an appropriate
simplex which contains a nondegenerate elliptic fixed point from the base
space. If instead of an elliptic fixed point, we have a stratum in
the base space which is closed and consists of nondegenerate elliptic singular points
of constant rank, then cutting the base space by an appropriate ``hyperplane" near this
stratum will lead to the symplectic blowing-up along the symplectic submanifold which
is the preimage of this stratum.

\begin{example}
Dehn surgery
\end{example}

Consider an integrable system with 2 degrees of freedom,
$\pi : (M^2,\omega, {\mathcal L}) \to O$, and let $D^2 \in O$
be a closed disk lying in the regular part of $O$. Cut out the piece
$\pi^{-1}(D^2)$ from the system, and then glue it back after some twisting.
This operation may  be called an {\it integrable Dehn surgery}, in analogy
with the well-known Dehn surgery in low-dimensional topology. It is
easy to see that any two roughly symplectically equivalent systems
with two degrees of freedom can be transformed to each other by performing
Dehn surgeries and then adding a magnetic term to the symplectic form.

\begin{example}
Hamiltonian Hopf bifurcation
\end{example}

There is a phenomenon called {\it Hamiltonian Hopf  bifurcation},
which happens in the Lagrange top and many other Hamiltonian systems
(see e.g. \cite{Meer}):  under some parameter change (e.g. the energy),
two pairs of purely imaginary eigenvalues (of the reduced linearized system
at an equilibrium) tend to each other, coincide, and then jump out of the
imaginary axe to form a quadruple of complex eigenvalues which is symmetric
with respect to both the real axe and the imaginary axe. One can verify,
for example, that a majority of integrable systems arising from the
so called {\it argument shift method} on coadjoint orbits of compact Lie
algebras (see e.g. \cite{MF}) exhibit such bifurcations. In terms of
integrable systems, this bifurcation amounts to a (generic) degenerate
corank 2 singularity which connects an elliptic-elliptic singularity
(i.e. a nondegenerate singularity of corank 2 which has 2 elliptic
components) to a focus-focus singularity in a 1-parameter family. From the
integrable surgery point of view, the passage from an elliptic-elliptic
singularity to a focus-focus singularity can be done via a small surgery, without
the need of an 1-parameter family. I will omit the explicit operation here.

\begin{example}
K3, ruled manifolds, etc.
\end{example}

It is easy to construct systems with 2 degrees of freedom for
which a topological 2-stratum $C$ of the base space is of any of
the allowed cases listed in Proposition \ref{prop:2domain}. The
most interesting case is $S^2$. The sphere $S^2$ admits an
integral affine structure with 24 singular points of focus-focus
type, which may be constructed as follows: Start with an integral
affine triangle (= base space of ${\Bbb C}P^2$ under a Hamiltonian
torus action). Cut from this triangle 3 small homothetic
triangles, each lying on one edge. Gluing together the pairs of
small edges that arise after we cut out the small triangles, we
obtain a triangle with an integral affine structure with 3
singular points of focus-focus type. We can glue 8 copies of this
new triangle together to obtain a sphere with 24 focus-focus
points. Proposition \ref{prop:realization2} shows that this $S^2$
is the base space of some integrable Hamiltonian system with 24
simple focus-focus singularities. Topologically, it is a torus
fibration over $S^2$ with 24 singular fibers of type $I^+$ in the
sense of Matsumoto, and the ambient space is diffeomorphic to a K3
surface (see \cite{Matsumoto} and references therein). We can also
go the other way around (less explicitly): start with a
holomorphic integrable system on a K3 surface (see
\cite{Markushevich}). Forgetting about the complex structure  and
taking the real part of the complex symplectic form, we get an
integrable system with 2 degrees of freedom whose base space is
homeomorphic to $S^2$. Similarly, we can construct a system whose
base space is homeomorphic to ${\Bbb R}P^2$, with 12 focus-focus
singular points. The ambient manifold will be diffeomorphic to an
Enriquez surface.

Assume now that the base space has no focus-focus singular point and is
diffeomorphic to the direct product of a graph or a circle with a closed interval
(the affine structure on $O$ needs not be a direct product). The ambient manifold of
integrable systems with such an orbit space $O$ are rational and ruled manifolds in
the sense of McDuff (see e.g. \cite{McDuff,MS}), which are analogs of complex ruled
surfaces (see e.g. \cite{BPV}). It can be shown easily that in this case, as in the
case of $S^2$ with 24 focus-focus points, we have $H^2(O,{\mathcal R}) = 0$ (for any
realizable affine structure on $O$). If $O$ is a product of 2 graphs which are not
trees, then there will be many topologically non-equivalent integrable systems over
$O$, because $H^2(O, {\mathcal R}) \neq 0$, etc.

\hspace{1cm}

{\bf Acknowledgements}. I would like to thank Jean-Paul Dufour,
Pierre Molino and L\^e H\^ong V\^an for conversations on
characteristic classes, Heinz Hansmann for explaining to me the
Hamiltonian Hopf bifurcation, Kaoru Ono for pointing out to me
Matsumoto's work on torus fibrations and the fact that an example
of mine of a symplectic 4-manifold is in fact diffeomorphic to a
K3 surface. I would like to thank also Michel Boileau, Lubomir
Gavrilov, Giuseppe Griffone and Jean-Claude Sikorav for an
invitation to Toulouse, for their very warm hospitality and
interesting discussions. And I'm especially indebted to Mich\`ele
Audin, Anatoly T. Fomenko and Alberto Verjovsky, without whom this
paper would not exist.

\bibliographystyle{amsplain}

\end{document}